\documentclass[a4paper,10pt,reqno]{amsart}
\usepackage[dvipsnames]{xcolor}
\usepackage[T1]{fontenc}
\usepackage{hyperref}
\hypersetup{
colorlinks = True,
linkcolor=Bittersweet,
citecolor=blue,
}
\definecolor{mycolor}{RGB}{226, 90, 104}
\definecolor{bluey}{rgb}{0.0, 0.0, 1.0}
\definecolor{greeny}{rgb}{0.0, 0.5, 0.0}
\definecolor{purpley}{rgb}{0.5, 0.0, 0.5}
\definecolor{orangey}{rgb}{1.0, 0.65, 0.0}
\definecolor{yellowy}{rgb}{1.0, 1.0, 0.0}
\definecolor{limegreeny}{rgb}{0.2, 0.8, 0.2}

\usepackage[nameinlink,capitalize, noabbrev]{cleveref}
\usepackage{multirow}
\usepackage{amsfonts,amsmath,mathrsfs,amssymb,amsthm,amscd,latexsym,amstext,amsxtra}
\usepackage{etoolbox}
\usepackage[all]{xy}
\usepackage{enumerate}
\usepackage{tabularx}
\xyoption{all}
\usepackage{srcltx} 
\usepackage[english]{babel}
\usepackage{mathtools}
\usepackage{epsfig}
\usepackage[]{epigraph}
\usepackage{xy}
\usepackage{enumerate}
\usepackage{graphicx}
\usepackage{float}

\DeclarePairedDelimiter{\floor}{\lfloor}{\rfloor}
\DeclarePairedDelimiter{\bracket}{\langle}{\rangle}

\newtheorem{theorem}{Theorem}[section]
\newtheorem{lemma}{Lemma}[section]
\newtheorem{corollary}{Corollary}[section]
\newtheorem{prop}{Proposition}[section]
\newtheorem{remark}{Remark}[section]
\newtheorem{Question}[theorem]{Question}
\newtheorem{definition}[theorem]{Definition}

\newtheorem{example}[theorem]{Example}

\DeclareMathOperator{\codim}{codim}

\newcommand{\nc}{\newcommand}
\nc{\cH}{{\mathcal H}}
\nc{\cA}{{\mathcal A}}
\nc{\cG}{{\mathcal G}}
\nc{\cC}{{\mathcal C}}
\nc{\cD}{{\mathcal D}}
\nc{\cO}{{\mathcal O}}
\nc{\cI}{{\mathcal I}}
\nc{\cB}{{\mathcal B}}
\nc{\cY}{{\mathcal Y}}
\nc{\cK}{{\mathcal K}}
\nc{\cX}{{\mathcal X}}
\nc{\cS}{{\mathcal S}}
\nc{\cE}{{\mathcal E}}
\nc{\cF}{{\mathcal F}}
\nc{\cZ}{{\mathcal Z}}
\nc{\cQ}{{\mathcal Q}}
\nc{\cN}{{\mathcal N}}
\nc{\cP}{{\mathcal P}}
\nc{\cL}{{\mathcal L}}
\nc{\cM}{{\mathcal M}}
\nc{\cT}{{\mathcal T}}
\nc{\cW}{{\mathcal W}}
\nc{\cU}{{\mathcal U}}
\nc{\cJ}{{\mathcal J}}
\nc{\cV}{{\mathcal V}}
\nc{\cR}{{\mathcal R}}
\nc{\bH}{{\mathbb H}}
\nc{\bA}{{\mathbb A}}
\nc{\bG}{{\mathbb G}}
\nc{\bC}{{\mathbb C}}
\nc{\bO}{{\mathbb O}}
\nc{\bI}{{\mathbb I}}
\nc{\bB}{{\mathbb B}}
\nc{\bY}{{\mathbb Y}}
\nc{\bK}{{\mathbb K}}
\nc{\bX}{{\mathbb X}}
\nc{\bS}{{\mathbb S}}
\nc{\bE}{{\mathbb E}}
\nc{\bF}{{\mathbb F}}
\nc{\bZ}{{\mathbb Z}}
\nc{\bQ}{{\mathbb Q}}
\nc{\bN}{{\mathbb N}}
\nc{\bP}{{\mathbb P}}
\nc{\bL}{{\mathbb L}}
\nc{\bM}{{\mathbb M}}
\nc{\bT}{{\mathbb T}}
\nc{\bW}{{\mathbb W}}
\nc{\bU}{{\mathbb U}}
\nc{\bD}{{\mathbb D}}
\nc{\bJ}{{\mathbb J}}
\nc{\bV}{{\mathbb V}}
\nc{\bbZ}{{\mathbb Z}}
\nc{\bR}{{\mathbb R}}
\nc{\fX}{{\mathfrak X}}
\nc{\Ida}{{\mathfrak a}}
\nc{\Idb}{{\mathfrak b}}
\nc{\Idp}{{\mathfrak p}}
\nc{\Idm}{{\mathfrak m}}
\nc{\fr}{{\rightarrow}}
\nc{\co}{{\nabla}}
\nc{\ok}{{\overline{K}}}
\nc{\Aut}{{\mbox{Aut}}}
\nc{\Int}{{\mbox{Int}}}
\nc{\la}{\longrightarrow}
\nc{\elemen}{\mathcal{E}(\bR^d)}
\nc{\Spec}{\mbox{Spec}}
\nc{\Proj}{\mbox{Proj}}
\nc{\Sym}{{\mathcal{S}}ym}
\nc{\divi}{\mbox{div}}
\nc{\Divi}{\mbox{Div}}
\nc{\Gal}{\mbox{Gal}}
\nc{\Bl}{\mbox{Bl}}
\nc{\GL}{\mbox{GL}}
\nc{\PGL}{\mbox{PGL}}
\nc{\mult}{\mbox{mult}}
\nc{\Conv}{\mbox{Conv}}
\nc{\Supp}{\mbox{Supp}}
\nc{\Pic}{\mbox{Pic}}
\nc{\cu}{{\overlineline{\nabla}}}
\nc{\con}{{\nabla}}
\nc{\vol}{{\mbox{vol}}}
\nc{\logc}{\bar{c}}

\patchcmd{\thebibliography}{\section*{\autorefname}}{}{}{}

\title{On the geography of $3$-folds via asymptotic behavior of invariants}

\author[Yerko Torres-Nova]{Yerko Torres-Nova}
\email{yerko.torres@usm.cl}
\address{Departamento de Matem\'atica, Universidad T\'ecnica Federico Santa Mar\'ia, Valpara\'iso, Chile.}

\begin{document}

\begin{abstract}
We study the geography problem for 3-folds of general type through the asymptotic behavior of invariants of $n$-th root covers. We first prove, in arbitrary dimension and for non-singular branch loci, that the Chern numbers are asymptotic to $n$ times the corresponding logarithmic Chern numbers of the base pair. In dimension three, for simple normal crossing branch divisors, we construct cyclic partial resolutions using toric methods and prove that, for asymptotic arrangements, the invariants $c_1^3, c_1c_2$, and $c_3$ have the same asymptotic behavior. We also obtain explicit families of 3-folds with ample canonical divisors that exhibit controlled Chern slopes.
\end{abstract}

\maketitle

\noindent \textbf{Keywords:} \textit{Complex algebraic geometry, Geography of threefolds, Chern numbers.}\\
\textbf{MSC2020:} \textit{14A99, 14C17, 14J30, 14J40.}

\tableofcontents

\section{Introduction}
We work with normal projective varieties $X$ over the complex numbers $\bC$. As usual, when $\dim X = 1$, $2$, or $d\geq 3$ we say that $X$ is a curve, a surface, or a $d$-fold respectively. The main purpose of this paper is to study the geography problem for $3$-folds of general type. Following the philosophy of Hirzebruch and Sommese \cite{hirze1983}\cite{somesse1984} studying the slopes of Chern numbers $c_1^2/c_2$ for minimal surfaces of general type, in \cite{hunt1989} the author proposed the study of slopes
$$[c_1^3,c_1c_2,c_3]\in\bP^2_{\bQ}$$
for non-singular minimal $3$-folds of general type. For a non-singular $3$-fold $X$, then we have Chern numbers $$ c_1^3(X) = -K_X^3,\quad c_1c_2(X) = 
24\chi(\cO_X),\quad c_3(X) = e(X).$$ 
For non-singular $3$-folds of general type, the triple $(K_X^3,\chi(\cO_X),e(X))\in \bQ\times \bZ^2$ is invariant among non-singular minimal models in its birational class \cite{matsuki2010} \cite{kollar1991}. In higher dimensions, minimal models may have singularities, and the usual Chern numbers are no longer directly available. Moreover, the topological Euler characteristic is not birationally invariant among singular minimal models. There are several inequalities for minimal $3$-folds, most of them for Gorenstein minimal $3$-folds, since we have well-behaved dualities on the cohomology groups. For example, in this case $\chi(\omega_X) = - \chi(\cO_X)$. The most important for us are the Miyaoka-Yau inequality  $-K_X^3 \geq 72\chi(\cO_X)$ (\cite{miyaoka87} \cite{Qyau2019}), and the classical general type inequalities $-K_X^3, \chi(\cO_X)< 0$, where the last one is just for the Gorenstein case. The search for Noether's type inequalities is an active research area, see \cite{chenhacon2006} \cite{yonghu2017} \cite{3noether2020}.\\

Previous results in this research line are the following. In \cite{hunt1989}, the author studied slopes of non-singular $3$-folds with ample canonical class in the chart $c_1c_2 \neq 0$. In \cite{chanlop} is proved that the Chern slopes $[c_1^3,c_1c_2,c_3]$ of non-singular $3$-folds with ample canonical bundle define a bounded region by finding inequalities for the coordinate $c_3/c_1c_2$. For a more general result, let $N(d)$ be the number of partitions of $d$. In \cite{dusun}, it is proved that Chern slopes $[c_1^d,\ldots,c_d]\in \bP^{N(d)-1}_{\bQ}$ of non-singular $d$-folds with ample canonical bundle define a bounded region in the chart $c_1^d \neq 0$. See \cite[p.14-18]{tesisYerko} to see an updated map in the chart $c_1c_2 \neq 0$ for non-singular $3$-folds with ample canonical class. The problem is that we do not know if those regions are optimal as we have for surfaces, i.e., the slopes $c_1^2/c_2$ are dense in the interval $[1/5,3]$. When working with slopes, it turns out that we can apply the tool of $n$-th root coverings to investigate their behavior. The main reason is that, as $n$ grows to infinity, the slopes erase the effects of the degree $n$ of the cover. This allows us to control the slopes from the base of the cover, and this is the main technique to get density results for surfaces.\\

Let $Z$ be a non-singular $d$-fold, and let $D_1,\ldots , D_r$ be distinct prime divisors in $Z$. Assume that $D_{red}  : = D_1+\ldots +D_r$ is a simple normal crossing divisor (SNC). For arithmetic simplicity, let $n>1$ be a prime number, and let $0<\nu_i<n$ be a collection of $r$ integers. Under certain conditions on the divisor $D = \nu_1D_1 + ... + \nu_rD_r$, we can construct an $n$-th root cover $h_n  \colon Y_n \to Z$ of degree $n$ \cite{viehwegesnault1992}. For curves, the prime divisors $D_1,\ldots ,D_r$ are distinct points on $Z$. Since points are isolated, we have $Y_n$ as a non-singular curve. By the Riemann-Hurwitz formula, we have $$\frac{c_1(Y_n)}{n} =  \frac{nc_1(Z) - (n-1)r}{n} = c_1(Z)-r+ \frac{r}{n}= \bar{c}_1(Z,D)+ \frac{r}{n}.$$ Here $\bar{c}_1(Z,D)$ is the first logarithmic Chern number of the pair $(Z,D_{red})$ (see \Cref{logprop}). Hence, if we fix the points $D_{red}$ and we consider partitions $\nu_1+\ldots +\nu_r=n$ with $n\gg 0$ a prime number, then we asymptotically have $c_1(Y_n) \approx n \bar{c}_1(Z,D)$. 

\begin{Question}
     Does this asymptotic phenomenon happen in higher dimensions?
\end{Question}

For $\dim Z \geq 2$, if the branch divisor $D_{red} = D_1+\dots+D_r$ has singularities, then $Y_n$ have rational singularities \cite{viehweg1977}. In order to have well-behaved invariants, we can choose a (partial) resolution of singularities. For $\dim Z = 2$, the asymptoticity of Chern numbers was proved in \cite{urzua1} for \emph{random}\index{random} $n$-th root covers. Let us explain briefly what random means (see  \Cref{asympres2} for more details). First, in dimension two, each singularity of $Y_n$ is a cyclic surface singularity of type $\frac{1}{n}(q,1)$ for some $0<q<n$. Thus, we use the Hirzebruch-Jung algorithm (see \Cref{HJ}) to resolve these singularities in a minimal way. For us, there are two important quantities: (1) the \emph{length of the resolution}, i.e., the number of steps of the algorithm, and (2) the Dedekind sums,
$$d(q,1,n) = \sum_{i=1}^{n-1} \left(\left( \frac{iq}{n}\right)\right) \left(\left( \frac{i}{n}\right)\right),$$
where $((\cdot)):\bR \to \bR$ is the saw-tooth function (see \Cref{dedekind}).
We get a resolution of singularities $X_n \to Y_n$ and the Chern numbers $c_1^2$ and $c_2$ depend on the lengths and the Dedekind sums coming from all cyclic singularities resolved. 
To guarantee asymptoticity, we have to consider
 \emph{asymptotic arrangements} defined in \Cref{asympres2}.\\
 
 The main set-up is when there exists an ample $H \in \Pic(Z)$ such that $D_{j} \sim H$ for each component of $D_{red}$. So, the $n$-th root cover construction applies when the $\nu_j's$ form a partition of $n$, i.e.,
 $$\nu_1+\ldots +\nu_r = n$$
 In \cite{urzua1} was proved that for a random partition $\nu_1+\ldots +\nu_r = n$, the probability of $D$ being an asymptotic arrangement tends to $1$ as $n$ grows.  In this way, for $n\gg 0$ we take random asymptotic partitions, and we get a family of random surfaces $X_n$ with
$$ c_{1}^2(X_n) \approx n\logc_1^2(Z,D),\quad c_2(X_n) \approx n\logc_2(Z,D).$$
In \cite{urzua2}, the asymptotic construction for surfaces was successfully completed in terms of the minimal model program \cite{mmp2010}. Explicitly, it was proved that the Chern numbers of the minimal model of each surface $X_n$ are asymptotically equal to the logarithmic Chern numbers of the base $(Z,D)$. Trying to do this in any dimension is the main goal of this research line. In this way, in higher dimensions, we have several issues with achieving an analog asymptotic result. For example, the singularities of $Y_n$ over triple intersections are not cyclic in general, and the choice of a right (partial) resolution of $Y_n$ with good behavior as $n$ grows is a challenging problem. 

\begin{Question}\label{question2}
Is there an analog of the asymptotic results in dimension two for dimension three?
\end{Question}

\subsection{Results of this work}

In \Cref{section0}, we prove that this phenomenon occurs in any dimension for logarithmic morphisms of degree $n$ with branch locus as a disjoint collection of non-singular distinct prime divisors $D_1,\ldots ,D_r$. In \Cref{pencilsec}, as an application, when the pair $(Z,D_{red})$ arises from a pencil through a complete intersection of hypersurfaces, we obtain a two-parameter family of accumulation points of Chern slopes of $3$-folds of general type with ample canonical divisor. The two limiting families lie on the lines
$$c_1^3 + c_3 = 2c_1c_2,\quad  c_1^3 + 3c_3 = 3c_1c_2,$$
and both converge to the point $[3,2,1]$.\\

We next turn to the SNC case of Question \ref{question2}. Here we first need to understand how the singularities of the cyclic covers and their resolutions affect the Chern numbers. The first result in this direction is in \Cref{asympchi}, where we find that the Chern number $c_1c_2 = 24\chi$ is \emph{asymptotic} and independent of the chosen resolution. The remaining Chern numbers require more care. The canonical self-intersection and the topological characteristic depend on the chosen resolution.  Since we desire at most an $O(n)$ behavior, the first issue is that the singularities of $Y_n$ over triple intersections have toric multiplicity $n^2$. This means that to connect with a non-singular model, at a bad choice of resolution we could have big exceptional data with respect to $n$, and so we would lose the asymptotic property. In \Cref{toricreso}, using toric methods, we construct a  \emph{local cyclic partial resolution}. In this way, we get singularities of multiplicity lower than $n$, and of cyclic quotient type. We are interested in cyclic quotient singularities since they are $\log$-terminal and have a well-known algorithm to resolve them: the Fujiki-Oka continued fraction (Cf. \cite{ashikaga2019}).  In \Cref{globalres} we globalize this local cyclic resolution, obtaining a cyclic partial resolution $X_n \to Y_n$. In \Cref{asympofK} and \Cref{asympofE} we show that the desired asymptotic property is satisfied. In terms of the minimal model program, it is important to have cyclic singularities, since they are well-understood. So, in future work we expect to discuss how these \emph{cyclic models} allow us to understand the minimal models.  We point out that the above is an embedded $\bQ$-resolution in the language of \cite{bart12}, introduced as an efficient resolution without \emph{useless data}.\\

\begin{small}
\noindent  \textbf{Acknowledgments:}  I am grateful to my Ph.D. thesis advisor Giancarlo Urz\'ua for his time,
guidance, and support throughout this work.  The results in this paper
are part of my Ph.D. thesis at the Pontificia Universidad Cat\'olica de
Chile.  I would also like to thank Jungkai Alfred Chen for his hospitality
during my stay at the National Center for Theoretical Sciences in Taiwan. Special thanks to
Pedro Montero and Maximiliano Leyton, for many comments and
suggestions to improve this work.  I was funded by the Agencia Nacional de Investigaci\'on
y Desarrollo (ANID) through the Beca Doctorado Nacional 2019.
\end{small}

\section{Preliminaries}

\subsection{Logarithmic properties}\label{logprop}
For a non-singular projective variety $Z$ of dimension $d$, let $A(Z) = \bigoplus_{e\geq 0} A^e(Z)$ be its Chow ring, where $A^e(Z)$ is the Chow group of cycles of codimension $e$. We denote its Chern classes by $c_e(Z) \in A^e(Z)$, i.e. the Chern classes of its tangent bundle $\cT_Z$ \cite[Appendix A]{hartshorne1977algebraic}. We have that $c_1(Z) = -K_Z$, where $K_Z$ is the canonical class. We define the Chern numbers of $Z$ as the degree of the top-intersection of its Chern classes
$$c_{i_1}(Z)\ldots  c_{i_m}(Z)\in A^d(Z), \quad i_1+\ldots +i_m = d.$$
In the following, if the context is understood, we abuse notation using the symbol $c_{i_1}\ldots  c_{i_m}(Z)$ for Chern numbers or more simply the notation $c_{i_1}\ldots  c_{i_m}$.\\

We have $c_1^d =  (-1)^dK_Z^d $, and since we work over $\bC$, it is well-known that $c_d = e(Z)$, the topological Euler characteristic. These numbers are codified into the Todd class of $\cT_Z$, and
as a consequence of the Hirzebruch-Riemann-Roch Theorem for $3$-folds we have the \emph{Noether's identity},
$$\chi(\cO_Z) = 
\dfrac{c_1c_2}{24}.$$

\begin{definition}
A \emph{simple normal crossing} (SNC) divisor $D = \sum_{j=1}^rD_j$ is a reduced effective divisor with distinct non-singular components $D_j$ satisfying the following condition: for each $p\in D$ there are local coordinates $x_1,\ldots ,x_d$ on $Z$  such that the equation defining $D$ on $p$ is $x_1\ldots x_e = 0$, with $e\leq d$.   
\end{definition}
 From \cite{iitaka1977} we introduce the following sheaf on $Z$.
\begin{definition}
For a SNC divisor $D$, the \emph{sheaf of $\log$-differentials along $D$}\index{sheaf of $\log$-differentials}, denoted by $\Omega^1_{Z}(\log D)$, is the $\cO_Z$-submodule of $\Omega^1_Z \otimes \cO_Z(D)$ described as follows. Let $p\in Z$ be a point.  
\begin{itemize}
    \item[(i)] If $p\not\in D$, then $ (\Omega^1_{Z}(\log D))_p = \Omega^1_{Z,p} $.
    \item[(ii)] If $p\in D$, we choose local coordinates $x_1,\ldots ,x_d$ on $Z$ with $x_1\ldots x_e = 0$ defining $D$ on $p$. Then, $(\Omega^1_{Z}(\log D))_p$ is generated as $\cO_{Z,p}$-module by 
$$\frac{dx_1}{x_1}, \ldots ,  \frac{dx_e}{x_e}, dx_{e+1},\ldots , dx_{d}.$$
\end{itemize}
If $D = \sum_{j=1}^r\nu_jD_j$ is a divisor on $Z$, whose associated reduced divisor  $D_{red}= \sum_j D_j$ is a SNC divisor, then for simplicity we set
$$\Omega_Z^1(\log D) := \Omega_Z^1(\log D_{red}).$$
\end{definition}

In the rest of this section, we assume $D = \sum_{j=1}^r\nu_jD_j$  as a divisor with $D_{red}$ a SNC divisor.

\begin{definition}
The \emph{$\log$-Chern classes}\index{$\log$-Chern classes} of a pair $(Z,D)$ are defined as
$$\logc_{i}(Z,D) = c_i(\Omega^1_{Z}(\log D)^{\vee}).$$
The \emph{$\log$-Chern numbers}\index{$\log$-Chern numbers} of a pair $(Z,D)$ are defined as the degree of top-dimensional intersections
$$\logc_{i_1}\ldots \logc_{i_m} := {\logc}_{i_1}(Z,D) \ldots {\logc}_{i_m}(Z,D)  ,\quad i_1+\ldots +i_m = d.$$
\end{definition}

 We set $\Omega^e_Z(\log D) := \bigwedge^e \Omega_{Z}^1(\log D)$ for any $1\leq e \leq d$. In this way, $\Omega^d_Z(\log D) = \cO_Z(K_Z + D_{red})$, i.e., ${\logc}_1(Z, D) = c_1(Z) - D_{red}$, and it is known that,
$${\logc}_d= e(Z) - e(D_{red}),$$
see \cite[Prop. 2]{iitaka1978}.

\begin{lemma}
We have a natural exact sequence
$$0 \to \Omega^1_{Z} \to  \Omega^1_{Z}(\log D) \to \bigoplus_{i=1}^r \cO_{D_i} \to 0,$$ which is known as the \emph{residual exact sequence}.
\end{lemma}
\proof{See \cite[Proposition 2.3]{viehwegesnault1992}. \qed}
\vspace{3mm}

From the residual exact sequence, we can compute the Chern polynomial through the identity
\begin{equation}\label{eq:polychern}
    c_t(\Omega_Z^1(\log D)) = c_t(\Omega_Z^1)\prod_{j=1}^r \left(\sum_{e=0}^d D_j^et^e\right).
\end{equation}
Let $1 \leq e \leq d$, and let $i_1+\ldots+i_m = e$ be any ordered composition of positive integers. We introduce the following notation,
$$D^{[i_1,...,i_m]} := \sum_{j_1<\ldots<j_m} D_{j_1}^{i_1}\ldots D_{j_m}^{i_m} ,\quad D^{[0]}=1.$$
Examples of this notation are
$$D^{[e]} = \sum_{j=1}^r D_j^e,\quad D^{[1,1]} = \sum_{j<k}D_jD_k,$$
and in general,
\begin{equation}\label{eq:identcomb}
    \prod_{j=1}^r\left(\sum_{e=0}^d D^e \right) =\sum_{e=0}^d \left( \sum_{i_1+\ldots+i_m = e} D^{[i_1,\ldots,i_m]}\right).
\end{equation}

A direct consequence of \Cref{eq:polychern} gives the following identity,
$$\logc_d(Z,D) = c_d(Z) + \sum_{e=1}^d (-1)^{e}c_{d-e}(Z)\left(\sum_{i_1+\ldots+i_m=e} D^{[i_1,\ldots,i_m]} \right).$$

\begin{lemma}\label{lemma0}
 Assume $D$ is non-singular, i.e., its non-singular components are pairwise disjoint.  Then for each $1\leq e\leq d$ we have $\logc_e(Z,D) = c_e(Z) + R_e(D)$ where 
 $$R_e(D)  = \sum_{\substack{k+l = e \\ k \neq e}} (-1)^lD^{[l]}c_k(Z),$$
 for each $e = 1,\ldots ,d$.  
\end{lemma}
\proof{Since $D$ is non-singular, then $D_iD_j = 0$ for all $i\neq j$. Thus, from the identity (\ref{eq:polychern}) , we get
$$\sum_{e=0}^d (-1)^e\logc_e(Z,D)t^e= \left(\sum_{e=0}^d (-1)^ec_e(Z)t^e\right)\left(\sum_{e=0}^d D^{[e]}t^e\right).$$
Using the Cauchy product formula for polynomials we have
$$(-1)^e\logc_e(Z,D) = \sum_{k+l = e} (-1)^{k}c_k(Z)D^{[l]},$$
and from this, the formula follows.
\qed}

\begin{lemma}\label{logchernum3}
Consider a non-singular $3$-fold $Z$, and $D = \sum_j \nu_jD_j$ on $Z$. We have
 $${\logc}_2(Z, D)  = c_2(Z) -D_{red}(c_1  - D_{red})  - D^{[1,1]}.$$
 Thus, the logarithmic Chern numbers for $3$-folds are,
 $$\logc_1^3 = c_1^3 -3 c_1^2 D_{red} + 3c_1([D^{[2]}+2D^{[1,1]}) - (D^{[3]} + 3D^{[1,2]} + 3D^{[2,1]} + 6D^{[1,1,1]})  $$
 $$\logc_1\logc_2  = c_1c_2 - D_{red}(c_1^2 + c_2) + c_1(2D^{[2]} + 3D^{[1,1]}) - D_{red}(D^{[2]} + D^{[1,1]}).$$
 $$\logc_3 =  c_3 - c_2D_{red} + c_1 (D^{[2]} + D^{[1,1]}) - \left(D^{[3]} + D^{[1,2]}+D^{[2,1]} + D^{[1,1,1]}\right).$$
\end{lemma}
\proof{ Using identity (\ref{eq:polychern}) for $d=3$, 
and looking for the degree $2$ terms we obtain $\logc_2(Z,D)$. The other formulas are direct computations using the above lemmas. 
\qed}

\begin{example}\label{exampleXgenemin}
Let $Z\hookrightarrow \bP^m$ be a non-singular projective $3$-fold. Let $H_1,\ldots,H_r \sim H$ hyperplane sections defining an SNC arrangement. We have
$$\logc_1^3 = c_1^3(Z) - 3rc_1^2(Z)H + 3r^2c_1(Z)H^2 -r^3\deg(Z)$$
$$\logc_1\logc_2 = c_1c_2(Z)-rH(c_1^2+c_2)(Z) + \left(2r + 3\binom{r}{2}\right)c_1(Z)H^2-\deg(Z)r\left(r + \binom{r}{2}\right),$$
$$\logc_3 = c_3(Z)-rHc_2(Z) + \left(r + \binom{r}{2}\right)c_1(Z)H^2 - \deg(Z)\left( r + 2\binom{r}{2} + \binom{r}{3}\right).$$
Thus as $r\to \infty$, we have
$$[\logc_1^3,\logc_1\logc_2,\logc_3] \longrightarrow [6,3,1]\in \bP_{\bQ}^2.$$
\end{example}

Logarithmic Chern classes are well-behaved under logarithmic morphisms.

\begin{definition}\label{logmorph}
Let $Z$ be a non-singular projective variety and $D$ an effective divisor with $D_{red}$ as SNC divisor. A surjective morphism $h\colon  Y \to Z$ between non-singular projective varieties is called a \emph{$\log$-morphism}, if $D'_{red} = (h^*D)_{red}$ is a SNC divisor. 
\end{definition} 

\begin{lemma}\label{lemmavieh}
 For any $\log$-morphism $h \colon (Y,D') \to (Z,D)$, we have an injection
$$h^*\Omega_{Z}(\log D) \hookrightarrow \Omega_{Y}(\log D').$$
Moreover, if $h$ is finite and ramified at $D$, then we have isomorphism outside the singularities of $D$.
\end{lemma}
\proof{See \cite[Lemma 1.6]{viehweg1982}. \qed}
\vspace{3mm}

\subsection{$n$-th root covers}\label{cyclicsec}
In this section, we follow \cite[Sec. 3]{viehwegesnault1992}.
Consider the following building data 
\begin{enumerate}[$\bullet$]
\item $Z$ is a non-singular projective variety of dimension $d$, 
\item $n\geq 2$ is a prime number,
\item $D = \sum_{j=1}^r\nu_j D_j$ with $0<\nu_j<n$ is an effective divisor on $Z$,
\item $D_{red} = \sum_{j=1}^rD_j$ is a SNC divisor, and
\item $D$ \emph{is divisible by $n$}, i.e., there is a line bundle $\cL$ on $Z$ such that $$\cO_Z(D) \simeq \cL^{\otimes n}.$$
\end{enumerate}
With this building data, we construct  
$$Y'_n = \Spec_Z \bigoplus_{i=0}^{n-1}\cL^{-i}$$
     with the cyclic algebra structure induced by the section defining $D$. The normalization  $f \colon Y_n \to Z$ of the natural morphism $Y_n' \to Z$ will be called the \emph{$n$-th root covering} associated to the building data. The morphism is finite of degree $n$, and $Y_n$ is projective and irreducible. The morphism $f$ is branched at $D$, and $Y_n$ has its singularities over the intersection of components of $D$. Thus, $Y_n$ is non-singular if $D_j\cap D_k = \emptyset$ for all $j,k$. On the other hand, if $D_{j_1}\cap\ldots\cap D_{j_e} \neq \emptyset$, then the branch divisor $D$ is defined by a local equation 
$$z_{j_1}^{\nu_{j_{1}}}\ldots z^{\nu_{j_{e}}}_{j_e} = 0,$$
where $z_1,\ldots ,z_d$ are local parameters for $Z$ on $p$. Thus the singularity of $Y_n$ over $D_{j_1}\cap\ldots\cap D_{j_e}$ is locally analytically isomorphic to the normalization of 
$$\Spec\left( \frac{\bC[z_1,\ldots ,z_d,t]}{t^n - z_{j_1}^{\nu_{j_{1}}}\ldots z^{\nu_{j_{e}}}_{j_e}} \right).$$ 

\begin{definition}
    A \emph{partial resolution of singularities}\index{partial resolution of singularities} of $Y_n$ is a projective, surjective, birational morphism $g:X \to Y_n$ with $X$ a projective normal variety having at most rational singularities. This last means that for any resolution of singularities $g' \colon X' \to X$ we have $R^ig'_*\cO_{X'} = 0$ for $i>0$. As usual, we omit the word \emph{partial} if $X$ is non-singular.
\end{definition}

Since the degree $n$ of $f$ is a prime number, we have $f^*(D_j) = n D'_j$, where $D_j' = f^*(D_j)_{red}$. Thus, for any partial resolution $h \colon X \to Y_n \to Z $ we must have a ramification formula
$$h^*(D_j) = nD'_j + \Delta_j,$$
where $\Delta_j$ is a divisor supported in the exceptional divisors of $h$ over $D_j$.

\begin{theorem}\label{maincyclic}
For the $n$-th root covering $f \colon Y_n \to Z$ we have
\begin{enumerate}
\item The morphism  $f$ is flat.
\item The variety $Y_n$ has rational singularities. 
\item We have the following decomposition on eigenspaces
$$f_*\cO_{Y_n} =  \bigoplus_{i=0}^{n-1}\cO_Z(-L^{(i)}),\quad L^{(i)} = iL - \sum_{j} \floor*{\frac{i\nu_j}{n}} D_j. $$
Indeed, $Y_n = \Spec f_*\cO_{Y_n}$. Therefore $f: Y_n \to Z$ is an affine morphism.
\item For any partial resolution of singularities $g \colon X \to Y_n$, the composition $h = f \circ g$ satisfies the following $\bQ$-numerical equivalence
$$K_X \sim_{\bQ} h^*\left( K_Z + \frac{n-1}{n}D_{red} \right) + \Delta,$$
where $\Delta$ is a divisor supported on the exceptional divisor of $g$.
\end{enumerate}
\end{theorem}
\proof{
For (1), (2), and (3) see  \cite[Ch. 3]{viehwegesnault1992}. For (4), we use the fact that $\codim(\mbox{Sing}(Y_n)) \geq 2$. In this way, there are non-singular open sets $U \hookrightarrow X \to V\hookrightarrow Z$ avoiding the singularities of $D$. From \Cref{lemmavieh}, we get $h^*\Omega^d_V(\log D) = \Omega^d_{U}(\log D')$ where $D'$ is the strict transform of $D$.  In terms of divisors, $$ h^*(K_V + D_{red}) = K_U + D'_{red} \sim_{\bQ} K_U + \frac{h^*(D_{red})}{n},$$ from where the result holds locally. After extending globally the exceptional term $\Delta$ appears.\\
\qed
}

\begin{lemma}\label{lemmaviewe}
    Let $Y$ be a normal variety and $g:X \to Y$ a proper, surjective, birational morphism. Assume that $X$ has rational singularities. Then $g_*\cO_X = \cO_Y$ and $R^ig_*\cO_X = 0$ for all $i>0$ if and only if $Y$ has rational singularities.
\end{lemma}
\proof{See \cite[Lemma 1]{viehweg1977}.\qed}

\begin{corollary}\label{constantchi}
For any partial resolution of singularities $g:X \to Y_n$, we have $$\chi(\cO_{X})=\sum_{i=0}^{n-1} \chi( \cO_Z(-L^{(i)})),$$ i.e., the analytic Euler characteristic of $X$ is independent of the chosen partial resolution.
\end{corollary}
\proof{ Since $Y_n$ has rational singularities, by \Cref{lemmaviewe} we have,
$$
g_*\mathcal O_X=\mathcal O_{Y_n},
\qquad
R^ig_*\mathcal O_X=0 \quad \text{for } i>0.
$$
Hence
$
\chi(\mathcal O_X)=\chi(\mathcal O_{Y_n}).
$
Since $f$ is finite,
$$
H^i(Y_n,\mathcal O_{Y_n})
\simeq H^i(Z,f_*\mathcal O_{Y_n}),
$$
and \Cref{maincyclic} gives
$$
\chi(\mathcal O_X)
=
\sum_{i=0}^{n-1}
\chi\bigl(\mathcal O_Z(-L^{(i)})\bigr).
$$\\
\qed
}
\subsection{Toric picture}\label{toricpic}

In this section, for toric varieties we mainly follow the notation of \cite{coxlittleschenk2011}.\\

Let $n>0$ be a prime number and $0\leq \nu_1,\ldots ,\nu_d<n$ integers. Choose a $\nu_k \neq 0$, and let $0\leq q_1,\ldots ,q_{d}<n$ be integers such that $\nu_j + q_j\nu_k \equiv 0$ modulo $n$. In particular, $q_{k}=n-1$. As usual set $N = \bZ^d$ and $N_{\bR} \cong \bR^d$ with canonical basis $e_1,\ldots ,e_d$, and $M = N^{\vee} \cong \bZ^d$ with $M_{\bR} = \bR^d$. Consider the semigroup 
$$S =\bracket*{e_1,\ldots ,e_{k-1},\sum_{j\neq k} q_je_j + ne_k ,e_{k+1},\ldots , e_d, \sum_{j\neq k} \frac{\nu_j + q_j \nu_k}{n}e_j + \nu_ke_k}_{\bN}. $$
Since $n$ is a prime number,
$$\sum_{j\neq k} \frac{\nu_j + q_j \nu_k}{n}e_j + \nu_ke_k = \sum_{j\neq k} \frac{\nu_j}{n} e_j + \frac{\nu_k}{n}\left( \sum_{j\neq k}  q_je_j +n e_k \right), $$
we have that the saturation \cite[p. 27]{coxlittleschenk2011} of $S$  is $S^{sat} = \sigma^{\vee} \cap \bZ^d$ where 
$$\sigma^{\vee} = C\left(e_1,\ldots ,e_{k-1},\sum_{j\neq k} q_je_j + ne_k ,e_{k+1},\ldots ,e_{d-1} , e_d\right)\subset M_{\bR}$$ is the simplicial $d$-cone defined by those elements. It is the dual cone of 
$$\sigma = C\left(ne_1 - q_1e_k , \ldots 
 ,ne_{k-1} - q_{k-1}e_k,e_k,ne_{k+1} - q_{k+1}e_k,\ldots  ,ne_{d} - q_{d}e_k \right)\subset N_{\bR}.$$
 Observe that $\mult(\sigma) = n^{d-1}$ and $\mult(\sigma^{\vee}) = n$.  Let $P_\sigma$ the fundamental parallelepiped of $\sigma$, i.e., the points of $\sigma$ with coordinates in $[0,1)$ respect its generators. Direct computations show that every element of $v\in P_{\sigma}\cap \bZ^d$ can written as
\begin{equation}\label{eq:decomv}
    v = \frac{\sum_{j\neq k}v_j(ne_j-q_je_k) + \left\{\sum_{j\neq k} v_jq_j \right\}_ne_k }{n}, \quad 0\leq v_j <n.
\end{equation}

\begin{prop}
The toric variety associated with the semigroup $S$ is
$$\Spec(\bC[S]) = \Spec\left( \frac{\bC[x_1,\ldots ,x_d,t]}{(t^n - x_1^{\nu_1}\ldots x_d^{\nu_d})}\right).$$
Moreover, its normalization corresponds to $\Spec(\bC[\sigma^{\vee} \cap \bZ^d])$.
\end{prop}
\proof{
For simplicity, we will prove the result in the case $k=d$. Since $Norm(\Spec(\bC[S])) = \Spec(\bC[S^{sat}])$ we will prove the first. Take the surjective morphism of semigroups $\phi \colon  \bN^{d+1} \mapsto S$ such that
$$ \phi(e_j) = e_j,\quad \phi(e_d) =\sum_{j=1}^{d-1} q_je_j +n e_d,\quad \phi(e_{d+1}) =\sum_{j=1}^{d-1} \frac{\nu_j + q_j \nu_d}{n}e_j + \nu_de_d.$$
It induces a surjective morphism of coordinate rings $f \colon  \bC[x_1,\ldots ,x_d,t] \to \bC[S]$, and by \cite{coxlittleschenk2011} in Proposition 1.1.9 it is known that 
$$\mbox{ker}(f) = (x_1^{a_1}\ldots x_d^{a_d}t^{a_{d+1}} = x_1^{b_1}\ldots x_d^{b_d}t^{b_{d+1}}  :  \phi(a) = \phi(b), a,b\in\bN^{d+1}).$$
If we set $x_j = a_j - b_j$, the condition $\phi(a) = \phi(b)$ gives equations
$$
\left\{\begin{array}{ccc}
x_1 + q_1x_d + \dfrac{\nu_1 + q_1\nu_d}{n} x_{d+1} &=& 0\\
\ldots  &\ldots &\ldots \\
x_{d-1} + q_{d-1}x_d + \dfrac{\nu_{d-1}+ q_{d-1}\nu_d}{n} x_{d+1} &=& 0\\
nx_d + \nu_d x_{d+1} &=&0
\end{array}
\right.$$
We can assume $x_{d+1} = nc$ with $c>0$, then $x_{d} = -\nu_d c$, and the equations reduce to
$$
\left\{\begin{array}{ccc}
x_1  &=& -\nu_1c\\
\ldots  &\ldots &\ldots \\
x_{d-1}  &=& -\nu_{d-1}c
\end{array}
\right.
\Rightarrow
\left\{\begin{array}{cccc}
b_j &=& a_j + \nu_j c,& 1\leq j\leq d  \\
a_{d+1} &=& b_{d+1} + nc&
\end{array}
\right.
$$
So $\mbox{ker}(f)$ is generated by elements of the form 
$$x_1^{a_1}\ldots x_d^{a_d}t^{b_{d+1}}((x_1^{\nu_1}\ldots x_d^{\nu_d})^c - (t^n)^c),$$
and the result follows.
\qed}

\begin{corollary}\label{coroindchi}
The normalization of the affine varieties $t^n = x_1^{\nu_1}\ldots x_d^{\nu_d}$ and $t^n = x_k\prod_{j\neq k}x_j^{n-q_j}$ are isomorphic.
\end{corollary}
\proof{ Observe that the cones defining both varieties are the same, equal to $\sigma$.
\qed}\\

\begin{remark} We point out the following.
    Assume that $k=e$ and $\nu_{e+1} = \ldots = \nu_{d} = 0$. Thus $q_{e+1} = \ldots = q_{d} = 0$, and we have a toric description of  the normalization of the varieties $t^n = x_1^{\nu_1}...x_e^{\nu_e}$  embedded in $\bA^d$ for any $1\leq e\leq d$.  
\end{remark}

\begin{remark}\label{remarkcyc}
It is known that the cone $\sigma^{\vee}$ defines a toric variety isomorphic to  the $d$-dimensional quotient cyclic singularity of type
$$\frac{(n-q_1,\ldots ,n-q_{k-1},1,n-q_{k+1},\ldots ,n-q_d)}{n}.$$
We have $\Spec\, \bC[\sigma \cap \bZ^d] \cong \bC^d/ \bracket{\phi}$, where $\phi \colon \bC^d \to \bC^d$ is defined by 
$$\phi \colon  (z_1,\ldots ,z_d) \mapsto (\zeta^{n-q_1}z_1,\ldots ,\zeta z_k,\ldots ,\zeta^{n-q_d}z_d),$$
with $\zeta^n = 1$ (Cf. \cite{ashikaga2015}). In this way, quotient cyclic singularities are geometrically dual to the singularities of $n$-th root covers. We will denote a cyclic singularity of this type by $C_{q_1,..., \hat{q_k},...,q_d}$. In dimension $2$, it occurs the accident that singularities of $n$-th root covers are also cyclic quotient singularities.
\end{remark}

\subsection{  Hirzebruch-Jung algorithm and Dedekind sums}
\subsubsection{Planar cones and Hirzebruch-Jung algorithm}\label{HJ}

Set $N\cong \bZ^d$ and $M = N^{\vee}$. A \emph{planar cone} $\tau$ in $N_{\bR}$ is a cone of dimension $2$, i.e., it is generated by two rays defined by primitive generators $v_0,v_{s+1}\in N$. Assume that $\mult(\tau) = n$. It is known that if $n>1$, then there exists some integer $0<q< n$ with $\gcd(n,q)=1$ and a
$$v = \frac{1}{n}(qv_0 + v_{s+1}) \in \tau \cap N$$
such that $C(v_0,v)$ is a non-singular cone. We say that $\tau$ is of \emph{type $\frac{1}{n}(q,1)$ in direction $v_0$ to $v_{s+1}$}. Consider the Hirzebruch-Jung algorithm of division for $n/q$, i.e., a pair of sequences
$$m_0 = n > m_1 = q > \ldots  > m_{s} = 1 > m_{s+1}=0,$$
$$n_0=0 <n_1 = 1 < \ldots  < n_{s} = q' < n_{s+1} = n ,$$
 which are related by 
$$m_{\alpha+1} = k_{\alpha}m_{\alpha}-m_{\alpha-1}$$
$$n_{\alpha+1} = k_{\alpha}n_{\alpha}-n_{\alpha-1},$$
where $k_1,\ldots ,k_s$ are integers satisfying $k_{\alpha} \geq 2$. Usually we denote $n/q= [k_1,\ldots ,k_s]$. These sequences define the Hirzebruch-Jung continued fraction as
$$\frac{n}{q}=k_1-\cfrac{1}{k_2-\cfrac{1}{\ddots - \frac{1}{k_s}}}.$$

\begin{lemma}\label{HJident}
For each $\alpha$, we have the following relations 
\begin{enumerate}
    \item $m_{\alpha}n_{\alpha+1} -m_{\alpha+1}n_{\alpha} =n, $
    \item $\mbox{gcd}(m_{\alpha},m_{\alpha+1}) = 1,$
    \item $\mbox{gcd}(m_{\alpha},n_{\alpha}) = 1.$
\end{enumerate}
\end{lemma}\qed

The non-singular resolution of the planar cone $\tau$ is a refinement by adding the rays defined recursively by
$$v_{\alpha} = \frac{m_{\alpha}v_{\alpha-1}+v_{s+1} }{m_{\alpha-1}} = \frac{m_{\alpha}v_0 + n_{\alpha}v_{s+1}}{n},\quad 1\leq \alpha\leq s.
 $$
Each cone $C(v_{\alpha},v_{\alpha+1})$ is non-singular, since 
 $$\det(v_{\alpha},v_{\alpha+1}) = \frac{1}{n}(m_{\alpha}n_{\alpha+1} -m_{\alpha+1}n_{\alpha}) = 1. $$

\begin{remark}\label{oppo}
Observe that sequence $n_{\alpha}$ is the sequence $m_{\alpha}$ for $n/q'$, i.e., if the pair $(m_{\alpha}',n_{\alpha}')$ is the resolution  of $n/q'$ then $(m_{\alpha}',n_{\alpha}')= (n_{s+1-\alpha},m_{s+1-\alpha})$. So we have a dual non-singular resolution given by the sequence
$$v_{\alpha}' = \frac{m'_{\alpha}v_{s+1}  + n'_{\alpha}v_0}{n},\quad 1\leq \alpha\leq s.$$
with $v_{\alpha} = v'_{s+1-\alpha}$.\\
\end{remark} 

After changing the $\mathbb Z$-basis of $N$ such that
$e_1=v$ and $e_2=v_0$, we have
$$
v_{s+1}=ne_1-qe_2.
$$
Hence, transversally to $\tau$, the associated toric variety is the cyclic
quotient surface singularity of type $\frac{1}{n}(q,1)$.
Therefore, the Hirzebruch-Jung subdivision above gives its minimal toric resolution.

\subsubsection{Dedekind sums}\label{dedekind} Consider the \emph{sawtooth function} $((x))\colon \bR \to \bR$ is defined as 
$$((x)) = \left\{ \begin{array}{cc}
x-\floor{x}-1/2 & x\in \bR\setminus \bZ\\
0 & x\in \bZ\end{array} \right. .$$ Observe that is an odd periodic function of period $1$. For $n\geq 3$ prime and $a_1,\ldots,a_d \in \bZ$ define the \emph{Dedekind sum} of dimension $d$ by
$$d(a_1,\ldots ,a_d,n) = \sum_{i=1}^{n-1} \left(\left( \frac{ia_1}{n}\right)\right) \cdots \left(\left( \frac{ia_d}{n}\right)\right).$$
By periodicity, we can reduce $a_i\in \bZ$ to $0\leq a_i < n$, then
$$d(a_1,\ldots ,a_d,n) =  d(\{a_1 \}_n,\ldots ,\{a_d \}_n ,n).$$
where $\{a_i \}_n$ is the residue modulo $n$ of $a_i$. Since $((x))$ is an odd function, we always have
$$d(-a_1,\ldots ,-a_d,n) =  (-1)^d d(a_1,\ldots ,a_d,n),$$
i.e., $d(a_1,...,a_d,n) = 0$ for any odd dimension $d$. We can rewrite this sum as
$$d(a_1,\ldots ,a_d,n) = \frac{1}{n^d}\sum_{i=1}^{n-1} \left( \{ ia_1\}_n - \frac{n}{2}\right) \cdots \left(\{ ia_d\}_n - \frac{n}{2}\right).$$

A well-known result due to Barkan \cite{Ba77} (see Holzapfel \cite{hol88}) is the following formula relating the length $s$ of the continued fraction and Dedekind sums,
\begin{theorem}\label{barkanhopls}
Let $n$ be a prime number, and $q$ be an integer such that $0<q<n$. Let $n/q = [k_1,\ldots ,k_s]$. Then
$$12d(1,q,n) + s = \sum_{\alpha=1}^s (k_{\alpha}-2) + \frac{q+q'}{n}.$$
\end{theorem}
\qed

\subsection{Asymptotic resolution in dimension 2}\label{asympres2}

Let $Z$ be a non-singular projective surface, and let $D$ be an effective divisor with SNC reduced divisor. Assume the necessary hypothesis to construct the normal $n$-th root cover $Y_n\to Z$ along $D$ (\Cref{cyclicsec}).  We have to choose a resolution of singularities $h \colon X_n\to Y_n \to Z$, and the Chern numbers $c_1^2,c_2$ of $X_n$ will depend on this resolution. The singularities of $Y_n$ over each point of $D_j\cap D_k$ are analytically isomorphic to the normalization of 
$$\Spec\left(\frac{\bC[x,y,t]}{t^n - x^{n-q_{jk}}y} \right),$$
where $\nu_j + q_{jk}\nu_k \equiv 0$ modulo $n$. This singularity is a cyclic quotient singularity of type  $\frac{1}{n}(q_{jk},1)$, and the singular point can be resolved by some weighted blow-ups. The exceptional data will be a chain of non-singular rational curves $\{E_1,\ldots , E_s \}$ with $E_jE_{j+1} = 1$ and $E_j^2 = -k_j$, where the $k_j\geq 2$ are the integers that define the negative regular continued fraction
$$\frac{n}{q_{jk}}=k_1-\cfrac{1}{k_2-\cfrac{1}{\ddots - \frac{1}{k_s}}},$$
usually called Hirzebruch-Jung continued fraction. See \Cref{HJ} for explicit computations. The number $s$ is called the \emph{length}\index{length} of the resolution and we denoted it by $\ell(q_{jk},n)$. In this way, we resolve all singularities of $Y_n$ obtaining a morphism $g \colon X_n\to Y_n$, with composition $h \colon  X_n \to Z$. In \cite{girstmair2003} and \cite{girstmair2006}, Girstmair proved that the lengths and the values of Dedekind sums have a particular asymptotical behavior.
 
\begin{theorem}[Girstmair]
\label{girsttheor}For $n\geq 17$ there exists a set $O_n \subset \{0,\ldots ,n \}$ such that for any $q\in O_n$ we have
$|d(1,q,n)| \leq 3\sqrt{n} + 5, \ell(q,n) \leq 3\sqrt{n} + 2$.
Moreover $|\{0,\ldots ,n\}\setminus O_n| \leq \sqrt{n}\log(4n)$.
\end{theorem}
\qed

\begin{remark}\label{asymptoticOk}
    Observe that from the Barkan-Holzapfel relation (\Cref{barkanhopls}), and combining it with the results of Girstmair, we have an asymptotic behavior for the coefficients of the Hirzebruch-Jung continued fraction in the following sense:  For a prime number $n\gg 0$, and integers $q\in O_n$ with $n/q = [k_1,\ldots ,k_s]$, we have
 $$\sum_{\alpha=1}^{s}(k_{\alpha}-2) = O(\sqrt{n}).$$
\end{remark}

\begin{definition}\label{asymtarrang}
A collection of prime divisors $\{ D_1,\ldots ,D_r\}$ on a non-singular $d$-fold $Z$ is an \emph{asymptotic arrangement}\index{asymptotic arrangements} if $D_{red} = D_1+\ldots+D_r$ is SNC, and for prime numbers $n\gg 0$:
\begin{enumerate}
    \item There are pairwise distinct multiplicities $0<\nu_j<n$, such that for any $j<k$ with $D_{j}\cap D_k \neq \emptyset$, we have $q_{jk} \in O_n$, the unique integer such that $\nu_j+q_{jk}\nu_k \equiv 0$ modulo $n$.
\item There are line bundles $\cL$ such that $\cO_Z\left( \sum_{j=1}^r \nu_jD_j \right) \simeq \cL^{\otimes n}$.
\end{enumerate}
\end{definition}

\begin{example}\label{partitions}
Inside the proof of \cite[Th. 6.1]{urzua1}, it was proved that for any large prime number $n$ there exists a partition $$\nu_1+\ldots +\nu_r = n,$$
with $q_{jk} \in O_n$ such that $\nu_j + q_{jk}\nu_k \equiv 0$ modulo $n$. We call it an \emph{asymptotic partition}\index{asymptotic partition} of $n$. Indeed, the probability of a partition of $n$ to be asymptotic tends to $1$ as $n$ grows to infinity. Thus, any collection of hyperplanes $\{H_1,\ldots ,H_r \}$ on $\bP^d$ defining a SNC divisor, is itself an asymptotic arrangement with
$$D = \nu_1H_1+\ldots +\nu_rH_r =(\nu_1+\ldots +\nu_r) H = nH,$$
 where $H$ is a general hyperplane section on $\bP^d$.
\end{example}
\vspace{1.5mm}

\section{Asymptoticity for a non-singular branch locus}\label{section0}
\subsection{Logarithmic asymptoticity}

\noindent Let $A(Z)_{\bR} = A(Z)\otimes_{\bZ} \bR$  be the extended Chow ring of a fixed non-singular variety $Z$ of dimension $d\geq 1$, and let $
D=D_1+\cdots+D_r
$
be a non-singular reduced divisor. For each $n\geq 1$ assume the existence of finite $\log$-morphisms $
h_n:(Y_n,D_n')\to (Z,D)
$ where 
$Y_n$ are non-singular varieties of the same dimension $d\geq 1$ with $\deg(h_n) = n$ and $D_n' = D_{n,1}'+\cdots+D_{n,r}'$ is the reduced inverse image of 
$D$.  We have a morphism of extended Chow rings $h_n^* \colon   A(Z)_{\bR} \to A(Y_n)_{\bR}$ for each $n$. Since $h_n$ is flat, we have that $h_n^*(A^e(Z)) \subset  A^e(Y_n)$ \cite[III.9.6]{hartshorne1977algebraic}, then the same applies for the extended ring. Assume moreover that for every $j$ there exists a positive integer
$e_{n,j}$ such that $
h_n^*(D_j)=e_{n,j}D_{n,j}'
$. Consider the branch divisor
$$
B_n:=\sum_{j=1}^r\frac{D_j}{e_{n,j}}
\in A^1(Z)_{\mathbb Q},
$$
satisfying $h_n^*(B_n) = D'_n$.

\begin{theorem}[$\log$-Ramification Formula]
 If
$
h_n^*\Omega_Z^1(\log D)
\simeq
\Omega_{Y_n}^1(\log D_n'),
$
then in $A^*(Y_n)_{\mathbb Q}$ we have the identity $$
c(Y_n)
=
h_n^*\left(
\overline c(Z,D)(1+B_n)
\right).
$$
Consequently, for every partition
$
i_1+\cdots+i_m=d,
$
one has the exact identity
\[
\frac{c_{i_1}\cdots c_{i_m}(Y_n)}{n}
=
\deg_Z\left[
\prod_{j=1}^m
\left(
\overline c_{i_j}(Z,D)
+
B_n\,
\overline c_{i_j-1}(Z,D)
\right)
\right].
\]
\end{theorem}
\proof{Since $D'_n$ is non-singular, its components are pairwise disjoint. 
From the residual exact sequence we have
$$
c(Y_n)
=
c\left(\Omega^1_{Y_n}(\log D'_n)^\vee\right)
\prod_{j=1}^r(1+D'_{n,j}).
$$
Since $D'_{n,j}D'_{n,k}=0$ for $j\neq k$, we get
$$
\prod_{j=1}^r(1+D'_{n,j})
=
1+D'_n.
$$
On the other hand, by hypothesis
$$
h_n^*\Omega^1_Z(\log D)
\simeq
\Omega^1_{Y_n}(\log D'_n),
$$
and therefore, by functoriality of Chern classes,
$$
c\left(\Omega^1_{Y_n}(\log D'_n)^\vee\right)
=
h_n^*\bar c(Z,D).
$$
Since $h_n^*(B_n)=D'_n$, we conclude that
$$
c(Y_n)
=
h_n^*\bar c(Z,D)(1+D'_n)
=
h_n^*\left(\bar c(Z,D)(1+B_n)\right).
$$
Looking at the component of degree $i$, we obtain
$$
c_i(Y_n)
=
h_n^*\left(
\bar c_i(Z,D)
+
B_n\bar c_{i-1}(Z,D)
\right),
$$
where $\bar c_{-1}(Z,D)=0$. Hence, for every partition
$i_1+\cdots+i_m=d$,
$$
c_{i_1}\cdots c_{i_m}(Y_n)
=
h_n^*
\left[
\prod_{j=1}^m
\left(
\bar c_{i_j}(Z,D)
+
B_n\bar c_{i_j-1}(Z,D)
\right)
\right].
$$
Finally, since $\deg(h_n)=n$, by the projection formula we get
$$
\frac{c_{i_1}\cdots c_{i_m}(Y_n)}{n}
=
\deg_Z
\left[
\prod_{j=1}^m
\left(
\bar c_{i_j}(Z,D)
+
B_n\bar c_{i_j-1}(Z,D)
\right)
\right].
$$\\
\qed
}

Now for each $n\geq 1$ assume there are parameters $0<\nu_j<n$ defining a divisor 
$$
D_n=\sum_{j=1}^r\nu_{n,j}D_j,
$$
divisible by $n$. Let
$
h_n:Y_n\to Z
$
be the non-singular $n$-th root covers associated to the building data $(Z,n,D_n)$. If $\nu_{n,j}$ and $n$ are coprime for all $j$, then $h_n^*(D_j) = n D_{n,j}'$, and we get the following asymptotic result.
\begin{corollary}
If $\gcd(\nu_{n,j},n)=1$ for every $j$, then for every partition
$i_1+\cdots+i_m=d$,
$$
\frac{c_{i_1}\cdots c_{i_m}(Y_n)}{n}
=
\deg_Z\left[
\prod_{j=1}^m
\left(
\overline c_{i_j}(Z,D)
+
\frac{D}{n}\,
\overline c_{i_j-1}(Z,D)
\right)
\right].
$$
In particular,
\[
\frac{c_{i_1}\cdots c_{i_m}(Y_n)}{n}
=
\overline c_{i_1}\cdots\overline c_{i_m}(Z,D)
+
O\left(\frac1n\right).
\]
\end{corollary}
\qed

For our purposes in geography, this result has a disadvantage: the difficulty in finding non-trivial pairs $(Z,D)$ whose covers $Y_n$ are minimal of general type. In \cite{bps2016} it was proved: Assume that $Z$ has a collection of disjoint divisors $\{D_j\}_{j\in J}$, if $|J|\gg 0$, then there exists a surjective morphism from $Z$ to a curve such that every $D_j$ is contained in a fiber. Thus, generically any variety $Z$ of $\dim Z \geq 2$ having collections of disjoint divisors is a fiber space $Z\to V$ over some variety $V$. 

\subsection{Geography Through Pencils}\label{pencilsec}

Following the previous discussion, let $V_0,V_{\infty} \subset \bP^d$ be hypersurfaces of degree $\delta\geq 1$ defined by sections $s_{0},s_{\infty} \in H^{0}(\bP^d, \cO_{\bP^d}(\delta))$. Consider the pencil 
$$\cP = \{V_{[\lambda_0,\lambda_{\infty}]} : [\lambda_0,\lambda_{\infty}] \in \bP^1 \},$$
where $V_{[\lambda_0,\lambda_{\infty}]}$ is the Cartier divisor defined by $\lambda_{\infty}s_0-\lambda_{0}s_{\infty} \in H^{0}(\bP^d, \cO_{\bP^d}(\delta))$. Let us assume that $V_0$ and $V_{\infty}$ intersect transversally, and that the base locus
$$B:= \bigcap_{[\lambda_0,\lambda_{\infty}] \in \bP^1 } V_{[\lambda_0,\lambda_{\infty}]} = V_{0}\cap V_{\infty},$$
is a non-singular and irreducible $(d-2)$-fold. So, we have a morphism
$$\phi : \bP^d\setminus B\to \bP^1,\quad x \mapsto [s_0(x),s_{\infty}(x)],$$
whose fibers are precisely the divisors of the pencil $\cP$. Let $Z=\mbox{Bl}_B \bP^d$ be the blowing-up through the base locus, so we get a birational morphism $\varphi : Z \to \bP^d $
inducing a fiber space
$\tilde{\phi}:  Z  \to \bP^1.$ Let $D = D_1 + \ldots+D_r$ be a non-singular reduced divisor whose reduced components $D_j$ are fibers of $\tilde{\phi}$, the strict transform of $V_j \in \cP$. In this situation, for each prime number $n>r$, choose positive integers
$0<\nu_1,\ldots,\nu_r<n$ such that
\[
\nu_1+\cdots+\nu_r=n.
\]
The divisors $D_j$ are fibers of the map $\tilde{\phi}$, so we assume that $F$ is a non-singular generic fiber and $D_j \sim F$. We have
\[
\sum_{j=1}^r \nu_jD_j \sim nF,
\]
and the previous theorem and corollary apply for the collection of cyclic covers
$$h_n : (Y_n, D_n') \to (Z, D).$$
In particular,
$$D \sim D_1 + \ldots + D_r \sim rF.$$
and, $F^e \sim 0$ for all $e>1$. Thus \Cref{lemma0} reduces to
$$\logc_e(Z,D) \sim c_e(Z) - rF\,c_{e-1}(Z),$$
and we have the following.
\begin{corollary}
    In the previous situation, we have
$$
\frac{c_{i_1}\cdots c_{i_m}(Y_n)}{n}
= c_{i_1}\cdots c_{i_m}(Z) - rF\sum_{j=1}^m c_{i_1}(Z)\ldots c_{i_j-1}(Z)\ldots c_{i_m}(Z)
+
O\left(\frac1n\right).
$$
\end{corollary}
\qed

To finish this section, let us give some explicit computations of our birational invariants chosen for $3$-folds. First, transversality implies that $V_j \sim \delta H$ has multiplicity $1$ along $B$ with $\varphi^*(V_j)  = D_j + E$. So, 
\begin{equation}\label{blowupeq}
    \delta \varphi^*(H) \sim  D_j + E \sim F + E.
\end{equation}
Moreover, it is known that $E^d = -(d-1)\delta^d$, so computing $(\delta \varphi^*(H))^d$ from \Cref{blowupeq} gives
$$E^{d-1}\,F = \delta^d,\quad \varphi^*(H)^{d-1} \,F  = \delta.$$

\begin{lemma}
    If $\delta> d+1$ and $n>r$, then the canonical divisor $K_{Y_n}$ is ample, i.e., $Y_n$ is a non-singular $d$-fold with ample canonical class. In particular,
    $$K_{Y_n}^d =  n \alpha^{d-1} (\alpha + d\delta  \beta ),$$
    where
    $$\alpha = \delta -d-1 ,\quad \beta = \left( \frac{r(n-1)}{n} - 1\right).$$
    
\end{lemma}
\proof{

Since the base locus $B$ is of codimension $2$, we have
$$K_Z = \varphi^*(K_{\bP^d}) + E = -(d+1)\varphi^*(H) + E,$$
where $H$ is an hyperplane section of $\bP^d$ and $E$ the exceptional divisor. So,
$$K_{Y_n} = h_n^*\left(K_Z + \frac{n-1}{n}D \right) \sim h_n^*\left(-(d+1)\varphi^*(H) + E + \frac{r(n-1)}{n} F\right)  \sim  h_n^*(\alpha  \varphi^*(H) +\beta F).$$
Note that both $\alpha$ and $\beta$ are positive rational numbers. Since $Z$ is the graph of the pencil, the morphism
$$
(\varphi,\tilde{\phi}): Z \hookrightarrow \mathbb{P}^d \times \mathbb{P}^1
$$
is a closed immersion. Hence,
\[
n(\alpha\phi^*H+\beta F)
=
(\varphi,\tilde{\phi})^*\mathcal{O}_{\mathbb{P}^d\times\mathbb{P}^1}(n\alpha,n\beta)
\]
is ample and also its pull-back by $h_n$. Finally, by the previous computations we have
\begin{align*}
    K_{Y_n}^d &= n(\alpha^d \varphi^*(H)^d + d \alpha^{d-1} \beta  \varphi^*(H)^{d-1} \,F) = n \alpha^{d-1} (\alpha + d\delta \beta ).
\end{align*}
\qed}

In what follows, denote the binomial polynomial as
$$\binom{d-\delta}{d} := \frac{1}{d!} \prod_{k = 1}^{d} (k - \delta).$$
\begin{lemma} The analytic and topological Euler characteristic of $Y_n$ are computed as
$$\chi(\cO_{Y_n}) = n + \frac{r(n-1)}{2} \left(\binom{d-\delta}{d}-1\right), $$
$$e(Y_n) = n(d+1 + e(B)) - (n-1)r \frac{(1-\delta)^{d+1} + \delta(d+1)-1}{\delta}.$$
\end{lemma}
\proof{ From \Cref{constantchi} we compute directly each $\chi( \cO_Z(-L^{(i)}))$ for $i=0,...,n-1$. In particular, we have
$$\cO_Z(-L^{(i)}) = \cO_Z(-r_iF),\quad r_i = \frac{1}{n} \sum_{j=1}^r \{ i\nu_j \}_n\in \bZ.$$
Since $F$ is a fiber of the map $\tilde{\phi} : Z \to \bP^1$, we have $\cO_Z(mF)|_F = \cO_F$ for any $m \in \bZ$. 
From the exact sequence
$$0 \to \cO_Z(-F) \to \cO_Z \to \cO_F \to 0,$$
we can compute inductively
$$\chi(\cO_Z(-r_iF)) = \chi(\cO_Z) - r_i \chi(\cO_F).$$
We have $\chi(\cO_{Z}) = \chi(\cO_{\bP^d}) = 1$, and 
$$\chi(\cO_{F}) =  \chi(\cO_{\bP^d}) - \chi(\cO_{\bP^d}(-\delta)) = 1 - \binom{d-\delta}{d}.$$
Since $i \mapsto \{ i\nu_j\}_n$ is a bijection of $\{1,...,n-1 \}$ when $\gcd(\nu_j,n) = 1$, we get
\begin{align*}
    \chi(\cO_{Y_n}) &= n - \left(1 - \binom{d-\delta}{d}\right)\sum_{i=1}^{n-1} r_i = n -\left(1 - \binom{d-\delta}{d}\right) \sum_{j=1}^{r} \frac{1}{n}\sum_{i=1}^{n-1}\{i\nu_j \}_n\\
    &= n -\left(1 - \binom{d-\delta}{d}\right) \frac{r(n-1)}{2}.
\end{align*}

The formula for $e(Y_n)$ comes from
\begin{align*}
    e(Y_n) &= e\left(Y_n \setminus \bigcup_{j=1}^r D'_j\right) + e\left(\bigcup_{j=1}^r D'_j\right) = ne\left(Z \setminus \bigcup_{j=1}^r D_j\right) + \sum_{j=1}^re(D_j)\\
    &= ne(Z) + (1-n)\sum_{j=1}^re(V_j)\\
    &= n(e(\bP^d) + e(B)) + (1-n)\sum_{j=1}^re(V_j).
\end{align*}
We have $e(\bP^d) = d+1$. Since all the $D_j$ are isomorphic to a degree $\delta$ hypersurface, we have
$$e(V_j)  = \frac{(1-\delta)^{d+1} + \delta(d+1)-1}{\delta},$$
from where the formula follows.
\qed
}

\begin{remark}
    Observe that in all cases, as $n\to \infty$ we get
$$\frac{K_{Y_n}^d}{n} \longrightarrow \alpha^d + (r-1) d\delta \alpha^{d-1},\qquad \frac{\chi(\cO_{Y_n})}{n} \longrightarrow 1 + \frac{r}{2} \left(\binom{d-\delta}{d}-1\right), $$
$$\frac{e(Y_n)}{n}  \longrightarrow d+1 + e(B) - r \frac{(1-\delta)^{d+1} + \delta(d+1)-1}{\delta}.$$
\end{remark}

In particular, for the case of $3$-folds we have that $B$ is a complete intersection curve of degrees $(\delta,\delta)$, so it is known that
$$K_B = 2(\delta - 2)H|_B,$$
with 
$$e(B) = -\deg K_B = 2(2-\delta)\delta^2.$$
Therefore, the above construction gives a two-parameter family of accumulation points.

\begin{corollary}  For every pair of integers $\delta>4$ and $r\geq 2$, the point
$$
P_{\delta,r}:=
\Big[
(\delta-4)^2\big((3r-2)\delta-4\big):
2r\delta(\delta^2-6\delta+11)-24:
r\delta(\delta^2-4\delta+6)+2\delta^3-4\delta^2-4
\Big]
\in \mathbb P^2_{\mathbb Q},$$
is an accumulation point of the Chern slopes $
[c_1^3,c_1c_2,c_3]$
of non-singular $3$-folds of general type with ample canonical class. 
\begin{enumerate}
    \item If $r$ is fixed and $\delta \to \infty$, then
$$
P_{\delta,r}\longrightarrow [\,3r-2:2r:r+2\,].
$$
on the line
$$c_1^3+ c_3=2c_1c_2.$$

\item For fixed $\delta$ and $r\to\infty$, we have
$$
P_{\delta,r}\longrightarrow
\Big[
3(\delta-4)^2:
2(\delta^2-6\delta+11):
\delta^2-4\delta+6
\Big],
$$
on the line
$$
c_1^3+3c_3=3c_1c_2.
$$
\end{enumerate}

Moreover, both families converge to the common point $[3,2,1]$
as $r\to\infty$ and $\delta\to\infty$, respectively.
\end{corollary}
\qed

\section{Asymptoticity of invariants}

\subsection{Asymptoticity of $\chi$ for $3$-folds}\label{asympchi}

Consider a data $(Z,D,n,\cL)$ as in \Cref{cyclicsec}, with $Z$ a non-singular projective $3$-fold. Let $h \colon  X_n \to Y_n \to Z$ be any resolution of singularities of the branched $n$-th root cover $Y_n$ along the effective divisor $D =\sum_{j=1}^r \nu_jD_j \sim nL$ whose reduced form is SNC. We have,
$$L^{(i)} = iL -\sum_{j=1}^r\floor*{\frac{i\nu_j}{n}} D_j = \frac{1}{n}\left(iD -\sum_{j=1}^rn \floor*{\frac{i\nu_j}{n}} D_j\right) = \frac{1}{n}\sum_{j=1}^r \{i\nu_j\}_n D_j.$$

\begin{lemma}\label{dederelations} Let $n$ be a prime number and let $0<a,b,c<n$. We have the following relations,
$$\sum_{i=1}^{n-1}  \{ia \}_n \{ib \}_n = n^2d(a,b,n) + \frac{n^2(n-1)}{4},$$
$$\sum_{i=1}^{n-1}  \{ia \}_n^2 \{ib \}_n = \sum_{i=1}^{n-1}  \{ia \}_n \{ib \}_n^2 = n^3d(a,b,n) + \frac{n^2(n-1)(2n-1)}{12} $$
$$\sum_{i=1}^{n-1}  \{ia \}_n \{ib \}_n\{ic \}_n = \frac{n^3}{2}\left( d(a,b,n) + d(a,c,n) + d(b,c,n) \right) + \frac{n^3(n-1)}{8},$$ 
\end{lemma}
\proof{We have an identity 
$$\sum_{i=1}^{n-1} \{ia \}_n^k = \sum_{i=1}^{n-1} i^k,$$
for each $k\geq 0$ integer. Then, we use repeatedly this identity in the following expressions. The first formula follows from,
$$n^2d(a,b,n) = \sum_{i=1}^{n-1} \left(\{ia\}_n - \frac{n}{2}\right)\left(\{ib\}_n - \frac{n}{2}\right),$$
and using that
$0 = \sum_{i=1}^{n-1} \left(\{ia\}_n - \frac{n}{2}\right)\left(\{ib\}_n - \frac{n}{2}\right)\left(\{ic\}_n - \frac{n}{2}\right),$
we get the other two.\\
\qed}

\begin{prop}
    We have, 
    $$\chi(\cO_{X_n}) = n\chi(\cO_Z) -\frac{1}{12}(R_1(n,D) + R_2(n,D) +R_3(n,D)),$$
    where
 $$R_1(n,D) =  \frac{(n-1)^2}{2n}D^{[3]} + \frac{(n-1)(2n-1)}{2n}(D^{[1,2]}+D^{[2,1]})+ \frac{3(n-1)}{2} D^{[1,1,1]},$$
$$R_2(n,D) = \frac{(1-n)}{2}c_1(Z)\left(\frac{(2n-1)}{n}  D^{[2]} + 3D^{[1,1]}\right) + \frac{(n-1)}{2}D_{red} (c_1^2(Z) + c_2(Z)),$$
\begin{align*}R_3(n,D) &=  6\left(\sum_{j<k}d(\nu_j,\nu_k,n)D_jD_k(D_j+ D_k+K_Z)\right.\\
&\left.+ \sum_{j<k<l} (d(\nu_j,\nu_k,n) + d(\nu_j,\nu_l,n) + d(\nu_k,\nu_l,n) )D_jD_kD_l\right).
 \end{align*}
\end{prop}
\proof{
By \Cref{constantchi}, and Hirzebruch-Riemann-Roch theorem for $3$-folds we can compute
\begin{align*}
\chi(\cO_{X_n}) &= \sum_{i=0}^{n-1}\chi( \cO_Z(-L^{(i)}))\\
&= n\chi(\cO_Z) - \frac{1}{12}\sum_{i=1}^{n-1}\left( L^{(i)}(L^{(i)} + K_Z)(2L^{(i)} +K_Z) + c_2(Z).L^{(i)} \right)\\
&=  n\chi(\cO_Z) - \frac{1}{12}\sum_{i=1}^{n-1}\left( 2(L^{(i)})^3 + 3(L^{(i)})^2K_Z + L^{(i)}K_Z^2 + c_2(Z).L^{(i)} \right)\\
&=n\chi(\cO_Z) -\frac{1}{12}R(n,D),
\end{align*}
where $R(n,D)$ is a quantity depending only on $n$ and $D$. We have the following identities: $\sum_{i=1}^{n-1} L^{(i)} = \frac{(n-1)}{2}D_{red}$, and $$\sum_{i=1}^{n-1} (L^{(i)})^2 = \frac{(n-1)(2n-1)}{6n}D^{[2]} + \frac{(n-1)}{2}D^{[1,1]} + 2\sum_{j<k}d(\nu_j,\nu_k,n) D_jD_k.
$$ The above is not difficult to deduce from the formulas in \Cref{dederelations}. To illustrate, we compute $\sum_i (L^{(i)})^3$ as follows. The first step, we compute explicitly:

\begin{align*}
(L^{(i)})^3 &= \frac{1}{n^3} \left( iD - \sum_{j=1}^r n\floor*{\frac{i\nu_j}{n}}D_j \right)^3\\
&=\frac{1}{n^3}\left(\sum_{j=1}^r \{i\nu_j \}_n D_j \right)^3\\
&= \frac{1}{n^3}\left( \sum_{j=1}^r \{i\nu_j \}_n^3 D_j^3 + 3\sum_{j<k} \{ i\nu_j \}_n^2\{ i\nu_k \}_n D_j^2D_k + \{ i\nu_j \}_n\{ i\nu_k \}_n^2 D_jD_k^2 \right.\\
& \hspace{12mm}+ \left. 6\sum_{j<k<l}\{ i\nu_j \}_n\{ i\nu_k \}_n\{ i\nu_l \}_n D_jD_kD_l \right)
\end{align*}
Applying directly the formulas in \Cref{dederelations}, we get.
\begin{align*}
\sum_{i=1}^{n-1} (L^{(i)})^3 &= \frac{(n-1)^2}{4n}D^{[3]} + \frac{(n-1)(2n-1)}{4n}(D^{[1,2]}+ D^{[2,1]})+ \frac{3(n-1)}{4}D^{[1,1,1]}\\
&+ 3\sum_{j<k}d(\nu_j,\nu_k,n)(D_j^2D_k+ D_jD_k^2)\\
& + 3\sum_{j<k<l} (d(\nu_j,\nu_k,n) + d(\nu_j,\nu_l,n) + d(\nu_k,\nu_l,n) )D_jD_kD_l.
\end{align*}
On the other hand, 

$$\sum_{i=1}^{n-1} (L^{(i)})^2K_Z = \frac{(1-n)}{2}c_1(Z)\left(\frac{(2n-1)}{3n}  D^{[2]} + D^{[1,1]}\right) + 2\sum_{j<k} d(\nu_j,\nu_k,n) D_jD_kK_Z,$$

$$ \sum_{i=1}^{n-1} L^{(i)}K_Z^2 =  \frac{(n-1)}{2}K_Z^2D_{red} =  \frac{(n-1)}{2}D_{red}c_1^2(Z), $$

$$ \sum_{i=1}^{n-1} L^{(i)}c_2(Z) = \frac{(n-1)}{2}D_{red} c_2(Z). $$
From here it is not too difficult to note that we have $R(n,D) = R_1(n,D) + R_2(n,D) +R_3(n,D)$, the quantities previously mentioned.\\
\qed}

\begin{theorem}\label{assympchi}
If $\{D_1,\ldots,D_r \}$ is an asymptotic arrangement, then
$$\frac{\chi(\cO_{X_n})}{n} \to \frac{\overline{c_1c_2}(Z,D)}{24},$$ 
as $n\to \infty$ for prime numbers $n\to \infty$.
\end{theorem}
\proof{ First observe the following limits
$$
\lim_{n\to \infty} \frac{R_1(n,D)}{n} = \frac{1}{2}D^{[3]} + (D^{[1,2]}+D^{[2,1]})+ \frac{3}{2} D^{[1,1,1]}= \frac{1}{2}D_{red}(D^{[2]} +D^{[1,1]}),
$$

$$
\lim_{n\to \infty} \frac{R_2(n,D)}{n} = -\frac{1}{2}c_1(Z)\left(2 D^{[2]} + 3D^{[1,1]}\right) + \frac{1}{2}D_{red} (c_1^2(Z) + c_2(Z)).
$$

\noindent From formulas in \Cref{logchernum3}, we get the identity
$$\lim_{n\to \infty} \frac{R_1(n,D)+R_2(n,D)}{n} =\frac{1}{2}( c_1c_2(Z) - \overline{c_1c_2}(Z,D)).$$
Since the collection of divisors is an asymptotic arrangement, we use \Cref{girsttheor} to get $R_3(n,D)/n\approx 0$ for prime numbers $n\gg 0$. Thus, in this case we have $$ \frac{\chi(\cO_{X_n})}{n} \to \frac{c_1c_2(Z)}{24} - \frac{1}{24}( c_1c_2(Z) - \overline{c_1c_2}(Z,D)) = \frac{\overline{c_1c_2}(Z,D)}{24},$$
for prime numbers $n\to \infty$.\qed}
\vspace{3mm}

\begin{remark}
Observe that the formula for $(L^{(i)})^3$ can be extended to higher dimensions. In the same way of \Cref{dederelations}, we can find formulas for $(L^{(i)})^e$ depending only on the combinatorial aspects of $D$ and higher dimensional Dedekind sums. Thus, for asymptoticity of $\chi(\cO_{X_n})$ in any dimension, we need asymptoticity of Dedekind sums, i.e., the higher dimensional analogs of Girstmair's results (\Cref{girsttheor}). For dimension $d\geq 4$ this is an open problem.   
\end{remark}

\subsection{Toric local resolutions}\label{toricreso}

\begin{figure}[h]
\begin{center}
\includegraphics[scale=0.4]{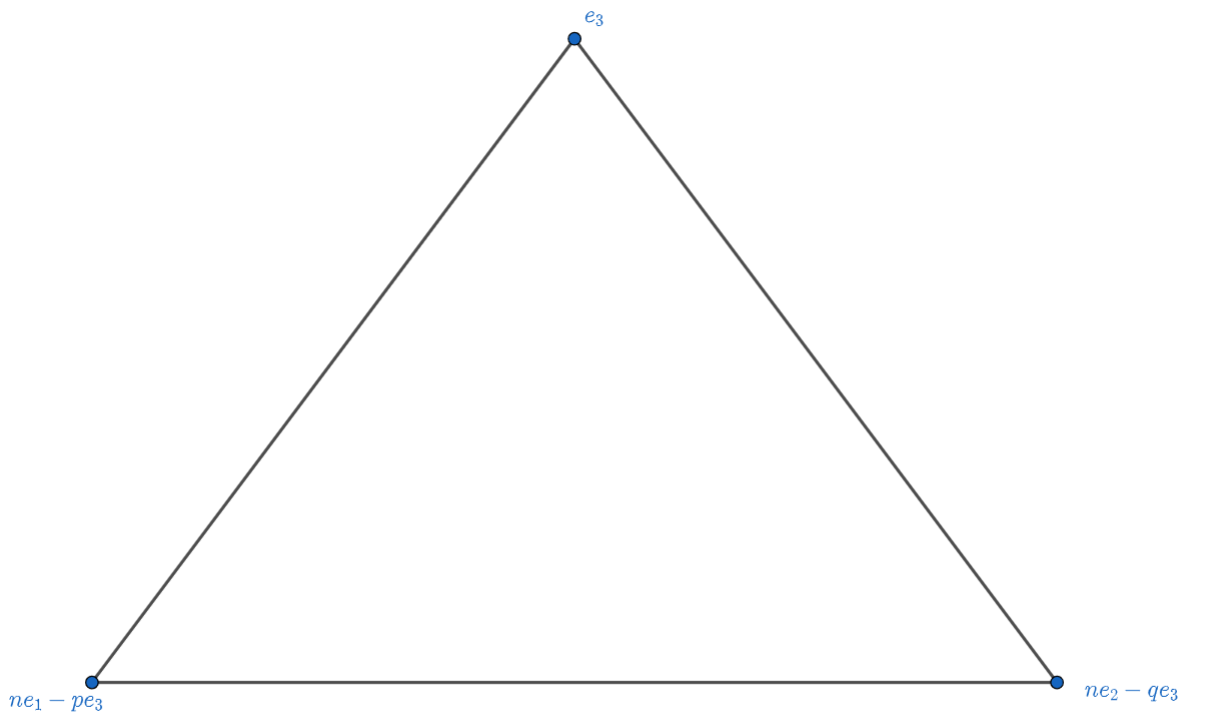}
\end{center}
\caption{Transversal section of $\sigma$.}
\label{fig:transcone}
\end{figure}

In this section, we study the $3$-fold singularity given by the normalization of $t^n = x_1^{\nu_1}x_2^{\nu_2}x_3^{\nu_3}$for $n\geq 2$ a prime number, and distinct integers $0<\nu_1,\nu_2,\nu_3<n$. The aim is to achieve good local resolutions of singularities in asymptotic terms with respect to $n$. Since this singularity is toric, it can be resolved by subdivisions of its associated cone obtaining a refinement fan. To assure asymptotic properties, we have to pay attention to the combinatorial aspects of the refinement. Let $0<p,q<n$ be integers such that $\nu_1+p\nu_3 \equiv 0$ and $\nu_2+q\nu_3 \equiv 0$ modulo $n$. By 
 \Cref{toricpic}, this $3$-fold singularity is a toric variety $Y_{p,q} := \Spec\bC[\sigma^{\vee}\cap \bZ^3]$ defined by the cone
$\sigma =  C(d_1,d_2,d_3)\subset \bR^3$, where $d_1 = ne_1-pe_3, d_2 = ne_2-qe_3$, and $d_3 = e_3$ are the primitive ray generators. A transversal section of this cone is sketched in \Cref{fig:transcone}. 

\begin{remark}
    This is the local model that will appear over a triple intersection
$D_j\cap D_k\cap D_l$, with
$p=q_{jl}$ and $q=q_{kl}$ after fixing $\nu_l$ as the distinguished
multiplicity. 
\end{remark}

\begin{figure}[h]
\begin{center}
\includegraphics[scale=0.4]{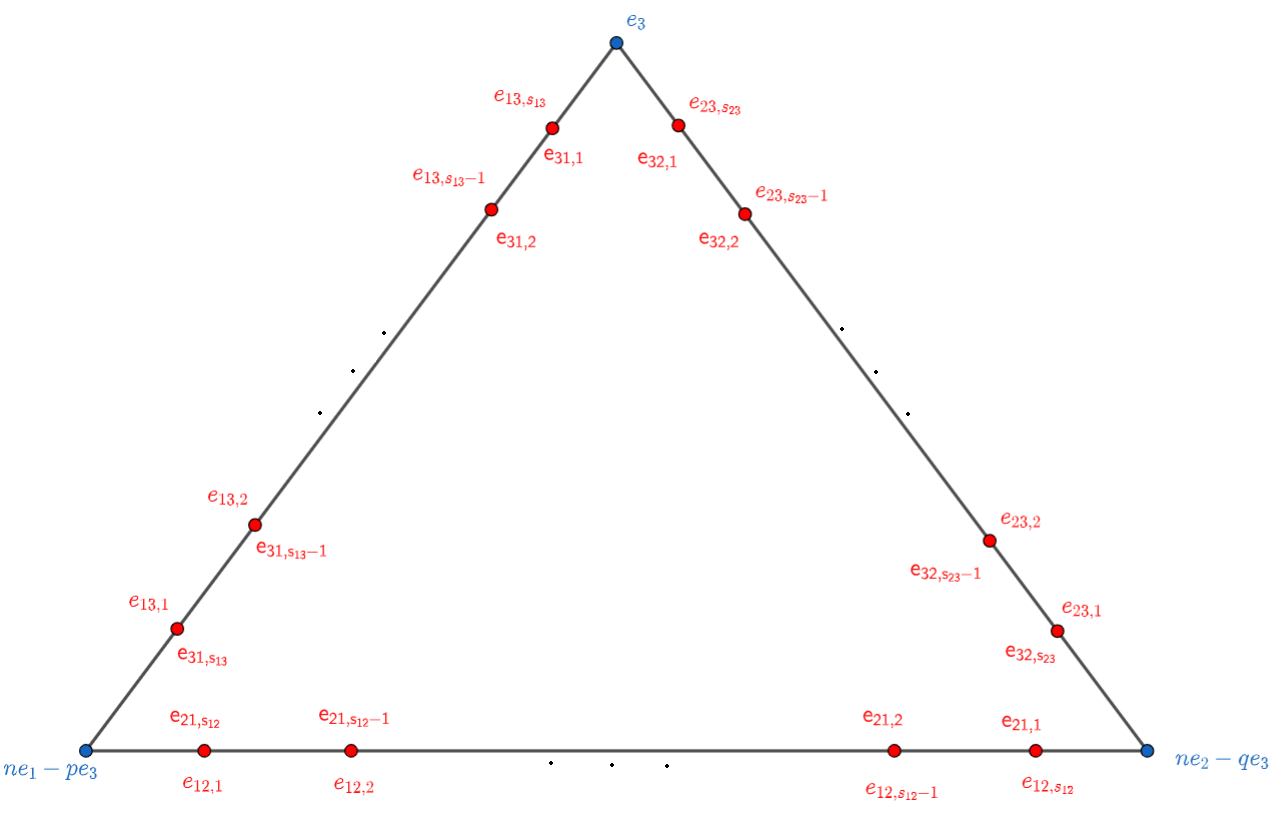}
\end{center}
\caption{Lattice points that resolve each wall}
\label{fig:resocone}
\end{figure}

The fan defined by the cone $\sigma$ has 2-dimensional faces (walls) given by $\tau_{jk}= C(d_j,d_k)$ for $j<k$. By \Cref{HJ}, each wall $\tau_{jk}$ can be resolved by Hirzebruch-Jung sequences  $(m_{jk,\alpha},n_{jk,\alpha})_{\alpha=0}^{s_{jk}+1}$ in direction $d_j$ to $d_k$ with initial data $$m_{jk,0} = n,\quad n_{jk,0} = 0,  \quad n_{jk,1}=1,\quad \forall j<k,$$
$$m_{13,1} = p',\quad m_{23,1} = q',\quad m_{12,1} = \{-p'q\}_n,$$
where $p'$ and $q'$ are the inverse modulo $n$ of $p$ and $q$. Thus, there are integers $k_{jk,\alpha} \geq 2$, such that
$$\frac{n}{m_{jk,1}} = [k_{jk,1},...,k_{jk,s_{jk}}].$$
Then, the walls $\tau_{jk}$ can be resolved by subdividing them in a sequence of steps by rays with generators  $(e_{jk,\alpha})_{\alpha=1}^{s_{jk}}$ defined recursively as
$$ e_{jk,\alpha} = \frac{m_{jk,\alpha}d_j+n_{jk,\alpha}d_k}{n},\quad 1\leq \alpha \leq s_{jk}.$$
 If we fix $j<k$, we denote by $e_{kj,\alpha}$ the ray generatorss in direction $d_k$ to $d_j$. We have the relation $e_{jk,\alpha} = e_{kj,s+1-\alpha}$. In \Cref{fig:resocone} are illustrated the border generators.

Canonical or minimal toric refinements do not provide a uniform combinatorial
model as $p,q$ vary. We therefore use an explicit refinement designed instead
to produce only cyclic quotient singularities while keeping the exceptional
data under asymptotic control. The advantage is that cyclic singularities have a well-known algorithm to resolve them \cite{fujiki74}\cite{ashikaga2019}.\\

From (\ref{eq:decomv}) we know that every $v\in P_{\sigma}\cap N$ can be written as
$$v = \frac{v_1d_1+v_2d_2+v_3d_3}{n}, \quad v_3 = \{pv_1+qv_2\}_n \quad 0 \leq v_i < n.$$ 
Fix a primitive lattice point $v\in P_\sigma\cap N$, with
$v_1+v_2+v_3\leq n$ and $v_i>0$.\\
 
\begin{center}
\begin{figure}[H]
\includegraphics[scale=0.4]{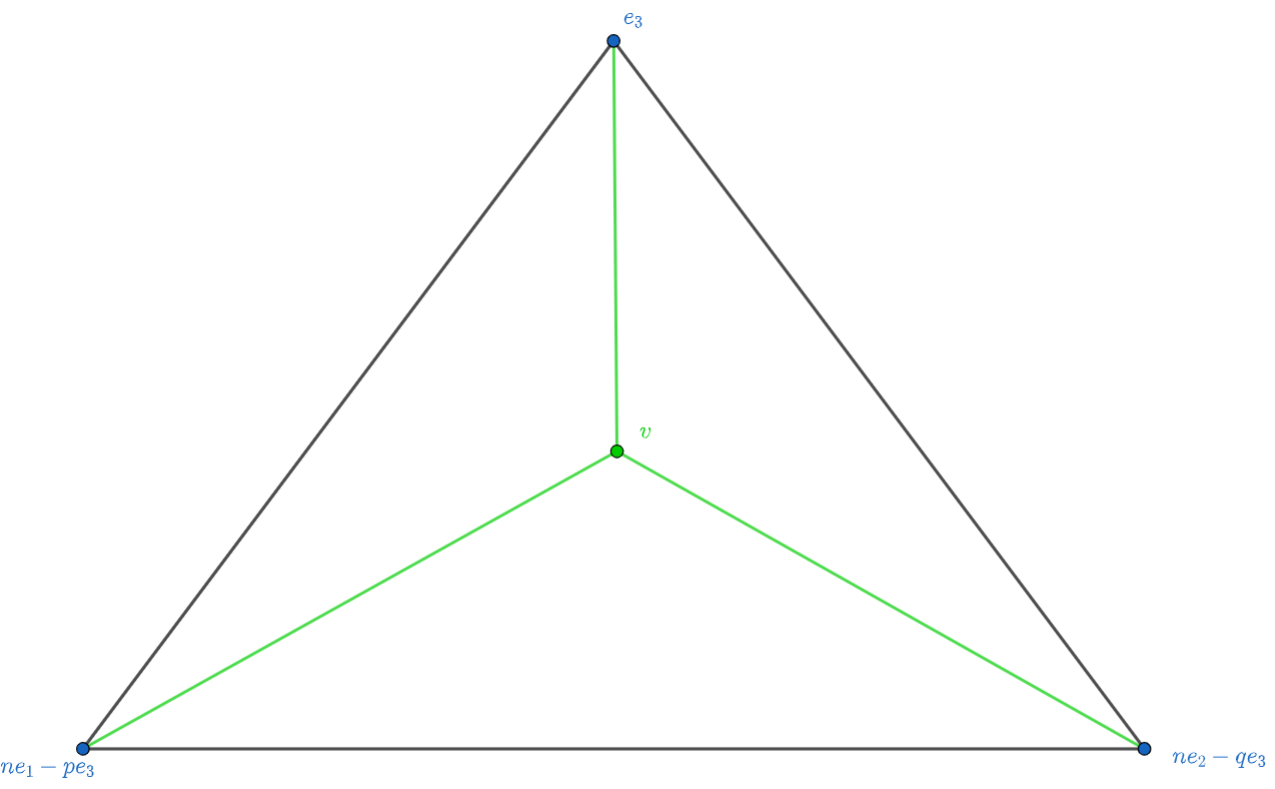}
\caption{Star subdivision along $v$}
\label{fig:centsub}
\end{figure}
\end{center}

\textbf{Step 1}: We refine by a star subdivision along $v$ obtaining a fan as illustrate \Cref{fig:centsub}.

\begin{figure}[h]
\begin{center}
\includegraphics[scale=0.4]{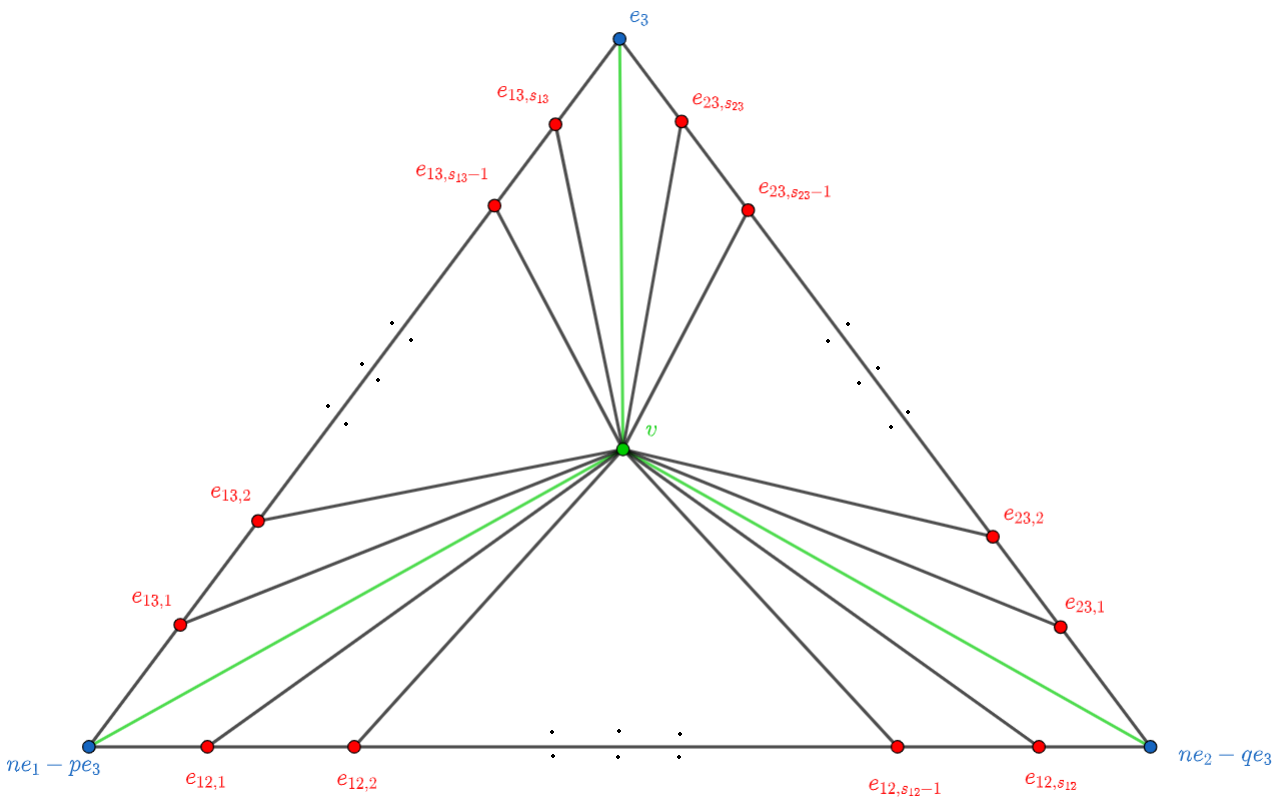}
\end{center}
\caption{Cyclic local resolution}
\label{fig:resoconecom}
\end{figure}

\textbf{Step 2}: We refine each wall by doing toric blow-ups following the Hirzebruch-Jung algorithm. So we obtain a refinement $\sigma^{*}$ and denote by $X$ the toric variety associated. This refinement gives us a birational projective morphism $g \colon  X \to Y_{p,q}$ \cite[11.1.6]{coxlittleschenk2011}. This refinement is sketched in \Cref{fig:resoconecom}.
Denote each $3$-cone of $\sigma^*$ by 
$$\sigma_{jk,\alpha} = C(v,e_{jk,\alpha},e_{jk, \alpha+1}),\quad 0 \leq \alpha \leq s_{jk}.$$
The $2$-cones of $\sigma^*$ are given in two types. The \emph{exterior walls} $\tau_{jk,\alpha} = C(e_{jk,\alpha},e_{jk,\alpha+1})$, and the \emph{inner walls} $\rho_{jk,\alpha} = C(v, e_{jk,\alpha})$. For any permutation $(v_j,v_k,v_l)$ with $j<k$, using determinants and properties of \Cref{HJ}, we have  
$$\mult(\sigma_{jk,\alpha})= v_l,\quad  \mult(\rho_{jk ,\alpha}) = \mbox{gcd}( v_jn_{jk,\alpha} - v_km_{jk,\alpha},v_l) ,\quad \mult(\tau_{jk,\alpha}) = 1.$$  
\begin{lemma}\label{cyclicsing}
Each cone $\sigma_{jk,\alpha}$ is a cyclic singularity of type
$$\frac{(a_{jk,\alpha}, b_{jk,\alpha},1)}{v_l},\quad a_{jk,\alpha} = \{m_{jk,\alpha+1}v_k - n_{jk,\alpha+1}v_j\}_{v_l},\quad b_{jk,\alpha}=\{m_{jk,\alpha}v_k - n_{jk,\alpha}v_j \}_{v_l},$$
where $\{\cdot \}_{v_l}$ is the residue modulo $v_l$.
\end{lemma}
\proof{ Since $\tau_{jk,\alpha}$ is non-singular, then $\sigma_{jk,\alpha}$ is semi-unimodular respect to $v$. By \cite[Prop. 1.2.3]{ashikaga2015} if there is positive integer $x,y$ such that
$$\frac{xe_{jk,\alpha} +ye_{jk,\alpha+1} + v}{v_l}\in \bZ^3,$$
then $x,y$ define the type of the cyclic singularity. Since $\gcd(n,v_l)=1$, we can solve the equations modulo $v_l$ and the result follows.\\ 
\qed
}

 Denote by $h\colon  X \to Y_{p,q} \to \bA^3$ the composition with the natural projection to $\bA^3$. Denote by $D_j$ the divisor in $\bA^3$ defined by the coordinate $x_j$. Each ray on $\sigma^*$ generated by $d_j,e_{jk,\alpha},$ or $v$ defines a toric divisor on $X$ given by $$\tilde{D}_j = V(C(d_j)),\quad E_{jk,\alpha} = V(C(e_{jk,\alpha})),\quad F =V( C(v)),$$ 
where $V(\cdot)$ denotes the orbit closure of a cone \cite[p. 121]{coxlittleschenk2011}. At the same time we fix notation for $j<k$ by $E_{kj,\alpha}  : = E_{jk,s+1-\alpha}$. 
 we can compute
$$h^*D_j = n\tilde{D}_j+ \sum_k\sum_{\alpha = 1}^{s_{jk}}m_{jk,\alpha}E_{jk,\alpha} + v_jF ,$$ 

\begin{equation}\label{eq:Klocal}
 K_X = g^*K_{Y_{p,q}} + \sum_{j<k}\sum_{\alpha = 1}^{s_{jk}}N_{jk,\alpha}E_{jk,\alpha}  + \frac{v_1+v_2+v_3-n}{n}F,   
\end{equation}
where 
$$N_{jk,\alpha}  : = \frac{m_{jk,\alpha} + n_{jk,\alpha} }{n} - 1.$$
It is satisfied the relation $k_{jk, \alpha}N_{jk, \alpha} - N_{jk, \alpha+1} = N_{jk, \alpha-1} - (k_{jk,\alpha}-2),$
which gives
$$\sum_{\alpha=1}^{s_{jk}} N_{jk,\alpha}(k_{jk,\alpha}-2) = -(N_{jk,1}+N_{kj,1}) - \sum_{\alpha=1}^{s_{jk}}(k_{jk,\alpha}-2).$$

\begin{prop}\label{propeff}
The $\bQ$-divisor  
$$h^*\left(\frac{n-1}{n}(D_1+D_2+D_3)\right) + \sum_{j<k}\sum_{\alpha}N_{jk,\alpha}E_{jk,\alpha} +\frac{v_1+v_2+v_3-n}{n}F,$$ is an  effective $\bZ$-divisor.
\end{prop}
\proof{The local pullback $h^*(D_1+D_2+D_3)$ equals
$$n(\tilde{D_1} + \tilde{D_2} + \tilde{D_3})+ \sum_{j<k,\alpha} (m_{jk,\alpha}+ n_{jk,\alpha}) E_{jk,\alpha} +(v_1+v_2+v_3)F.$$
Thus, $h^*\left(\frac{n-1}{n}(D_j+D_k+D_l)\right) + \Delta$ equals to
$$(n-1)(\tilde{D_1} + \tilde{D_2}+\tilde{D_3})+ \sum_{j<k,\alpha} \left(m_{jk,\alpha}+ n_{jk,\alpha}-1\right) E_{jk,\alpha} + (v_1+v_2+v_3-1)F,$$
i.e., an effective $\bZ$-divisor. 
\qed
}
\vspace{3mm}

 The following arithmetic lemma will be useful in \Cref{globalres}. 
\begin{lemma}\label{sumnasym}
Assume that $p\neq q$ and that $p,q\neq n-1$.
Then there exists a primitive $v\in P_\sigma\cap N$ such that
\[
v_1+v_2+v_3=n.
\]
In particular, the coefficient of $F$ in \eqref{eq:Klocal} is zero.
Moreover, for every fixed $C>4$ and $n\gg0$, the point $v$ can be
chosen so that
\[
\frac{v_j}{v_k}\le C
\]
for every $j,k$.
\end{lemma}
\proof{
Set
$
c=\{-(p+1)(q+1)'\}_n.
$
Since $p,q\neq n-1$ and $p\neq q$, we have
$c\not\equiv 0,-1\pmod n$. For $1\le a<n$, define
\[
v_1(a)=a,\qquad
v_2(a)=\{ca\}_n,\qquad
v_3(a)=\{-(1+c)a\}_n.
\]
Since
$
\{(p+1)+(q+1)c\}_n = 0,
$
we have
\[
pv_1(a)+qv_2(a)\equiv v_3(a)\pmod n,
\]
so $v(a)\in P_\sigma\cap N$. Moreover,
\[
v_1(a)+v_2(a)+v_3(a)\equiv0\pmod n.
\]
Since all three coordinates belong to $\{1,\ldots,n-1\}$, their
sum is either $n$ or $2n$. On the other hand,
\[
v_i(n-a)=n-v_i(a),
\]
which implies
\[
\sum_i v_i(n-a)=3n-\sum_i v_i(a).
\]
Therefore exactly one element in each pair $\{a,n-a\}$ satisfies
\[
v_1(a)+v_2(a)+v_3(a)=n.
\]
Thus there are $(n-1)/2$ such lattice points. Fix $C>4$ and set $\epsilon=(C+2)^{-1}$. For each $i$, the map
$a\mapsto v_i(a)$ permutes $\{1,\ldots,n-1\}$. Hence fewer than
$\epsilon n$ of the above values satisfy $v_i(a)<\epsilon n$.
Consequently fewer than $3\epsilon n$ values have at least one
coordinate smaller than $\epsilon n$. Since
\[
3\epsilon<\frac12,
\]
for $n\gg0$ we can choose one of the $(n-1)/2$ points satisfying $
v_i\ge\epsilon n,\, i=1,2,3.$ Since their sum is $n$, also
\[
v_i\le (1-2\epsilon)n=C\epsilon n.
\]
Therefore
\[
\frac{v_j}{v_k}\le C
\]
for every $j,k$. Finally, if $v=mw$ for
some $w\in N$ and $m>1$, then, since $w\in P_\sigma$, uniqueness of
the coordinates in $P_\sigma$ implies that $m$ divides each $v_i$.
Hence $m$ divides $v_1+v_2+v_3=n$. Since $n$ is prime and
$0<v_i<n$, this is impossible unless $m=1$. Thus, $v$ is primitive.\\
\qed}

\subsection{Local Intersection Formulas}

Each inner wall defines a closed curve on $X$ given by 
$$C_j = V(C(d_j,v)),\quad C_{jk,\alpha} = V(C(\rho_{jk,\alpha})).$$
The refinement $\sigma^*$ is simplicial, hence $X$ is $\bQ$-factorial and $K_{X}$ is $\bQ$-Cartier.  For any pair $v_j,v_k$, $j<k$, let $v_l$ be the another coordinate.  The following relation among lattices generators,
$$e_{jk,\alpha-1} + (-k_{jk,\alpha})e_{jk,\alpha} + 0\cdot v +e_{jk,\alpha+1} = 0,$$
$$v_j e_{lj,1} + \left(\frac{v_l-m_{lj,1}v_j-m_{lk,1}v_k}{n}\right)d_l +  (-1)v + v_ke_{lk,1}=0,$$   describe the intersection theory on $X$ \cite[Ch. V]{fulton1993introduction}. We have,
$$E_{jk,\alpha}C_{jk,\alpha \pm 1} = \frac{\mult(\rho_{jk,\alpha})}{v_l} ,\quad E_{jk,\alpha}C_{jk,\alpha} = -\frac{k_{jk,\alpha}\mult(\rho_{jk,\alpha})}{v_l},\quad FC_{jk,\alpha} = 0,$$
$$\tilde{D}_l\, C_l = \frac{\gcd(v_j,v_k)(v_l - m_{lj,1}v_j -  m_{lk,1}v_k) )}{nv_jv_k},\quad  F\, C_l =-\frac{\gcd(v_j,v_k)}{v_jv_k} ,$$
$$K_X\, C_l = -\frac{\gcd(v_j,v_k)}{v_jv_k}\left(v_j+v_k-1+\frac{v_l-m_{lj,1} v_j - m_{lk,1}v_k}{n} \right), $$
\begin{equation}\label{eq:Kcurvasjk}
   K_X\, C_{jk,\alpha} =  \frac{\mult(\rho_{jk,\alpha})}{v_l}(k_{jk,\alpha}-2),\quad 1\leq \alpha \leq s_{jk}. 
\end{equation}

\begin{lemma} We have
$$F^3 = \frac{n}{v_1v_2v_3}.$$
\end{lemma}
\proof{The divisor $v_jF$ is Cartier for any $j$, then from the pullback identities above we have
$$
v_1v_2v_3F^3 = \prod_{j=1}^3 \left(h^*D_j - n\tilde{D}_j - \sum_k\sum_{\alpha} m_{jk,\alpha}E_{jk,\alpha}\right) = h^*D_1\, h^*D_2\, h^*D_3 = n,$$
where the last is by projection formula.\qed
}
\vspace{3mm}

As a consequence, we can compute,
\begin{align*}
K_XF^2 &= -F^2\left(D_1+D_2+D_3 + F\right)\\
&= \frac{v_1+v_2+v_3 - n}{v_1v_2v_3}.
\end{align*}
From the last one, we get
\begin{align*}
K_X^2\, F &=-\sum_l \frac{K_X\, C_l}{\gcd(v_j,v_k)} - \sum_{j<k,\alpha} \frac{K_X\, C_{jk,\alpha}}{\mult(\rho_{jk,\alpha})} - K_XF^2.
\end{align*}

\begin{example}\label{examplepq1}
\emph{Case $\{p + q\}_n = 1$:}  For $v_1 = v_2 = 1$, we have $v_3 = \{ p+q\}_n =1$. Thus, these coordinates define an interior lattice point $v$, and it minimizes $v_1+v_2+v_3$. We have $$\mult(\sigma_{jk,\alpha})= \mult(\rho_{jk ,\alpha}) =  \mult(\tau_{jk,\alpha}) = 1,$$
i.e., $X$ is non-singular. On the other hand,
$$K_XC_{jk,\alpha} = (k_{jk,\alpha} - 2),\quad K_XC_l = 0,$$
for all $j,k$ and $l$. Thus, $K_X$ is nef.
\end{example}

\subsection{Global resolution}\label{globalres}

 Let $\{D_1,...,D_r\}$ be an asymptotic arrangement on $Z$.  Thus, for prime numbers $n\gg 0$ we have multiplicities $0<\nu_j<n$ depending on $n$, with its respective $q_{jk}\in O_n$. We have $n$-th root covers $f_n \colon  Y_n \to Z$ branched at each $D = \sum_j\nu_j D_j$. Let us fix a $n\gg 0$, so we drop the subscript $n$ of the morphisms, i.e., we are fixing a $f\colon  Y_n \to Z$.\\
 
 For $j<k$, the singularities of $Y_n$ over curves in $D_{jk} := D_j\cap D_k$ are locally analytically isomorphic to $C_{q_{jk},n}\times \bC$ where $C_{q_{jk},n}$ is the surface cyclic quotient singularity of type $\frac{1}{n}(q_{jk},1)$. For a triple $j<k<l$, the singularity of $Y_n$ over a point in $D_{jkl}:= D_j\cap D_k\cap D_l$ is locally analytically isomorphic to the normalization of $\Spec\bC[\sigma_{jkl}^{\vee}\cap M ]$ where $\sigma_{jkl}$ is a cone with walls of types $\frac{1}{n}(q_{jk},1), \frac{1}{n}(q_{jl},1), \frac{1}{n}(q_{kl},1)$ as we see in \Cref{toricreso}. We will get the \emph{cyclic resolution}\index{cyclic resolution} $X_n \to Y_n$ via weighted blow-ups in two steps.\\

 \textbf{Step 1:} For a triple intersection $D_{jkl}\neq\varnothing$,
notice that the hypotheses of \Cref{sumnasym} are satisfied for $n\gg 0$. Indeed, $p,q\neq n-1$ and cannot be equal since the $\nu_j's$ are pairwise distinct. Thus, for each triple intersection $D_{jkl}$, we choose the interior
lattice point $v^{jkl}$ provided by the lemma. Since singularities over $D_{jkl}$ are isolated, for each $\sigma_{jkl}$, we do a weighted blow-up at $v^{jkl}$. So, this refinement locally gives a projective morphism which is a blow-up over an isolated point \cite[11.1.6]{coxlittleschenk2011}. In this case, we get a projective birational morphism $h':X_n' \to Y_n$.  We have exceptional divisors $F_{jkl}$ whose components are over the points of $D_{jkl}$ and they are isomorphic to weighted fake projective planes \cite{buczynska2008fake}. For future computations, we fix the notation of $v^{jkl}$ and $F_{jkl}$ independent of the order of the triple $j,k,l$. For example, $v^{jkl} = v^{kjl}$.\\
 
 \textbf{Step 2:} Away from the triple intersections, the singularities over $D_{jk}$ are unchanged. These are locally analytically isomorphic to
$C_{q_{jk},n}\times\mathbb C$. We resolve them by the Hirzebruch-Jung
subdivision associated with $q_{jk}$. Near a triple intersection, this subdivision agrees with the subdivision of
the corresponding boundary wall of the local toric model constructed in
\Cref{toricreso}. Hence these local refinements are compatible on overlaps and glue
to a global projective birational morphism. Since $Y_n$ is irreducible and the construction consists of successive
toric weighted blow-ups, $X_n$ is irreducible.
 
\begin{prop}[and Definition] There exists a \emph{cyclic partial resolution} $g:X_n \to Y_n$, i.e. a projective, surjective, birational morphism such that $X_n$ is irreducible and, it has at most cyclic quotient singularities of order lower than $n$. 
\end{prop}
\qed

By construction, the local fans are simplicial. Hence $X_n$ is
$\mathbb Q$-factorial and $K_{X_n}$ is $\mathbb Q$-Cartier. In the rest, we will abuse notation using $D_{jk}$ and $D_{jkl}$ for both, set-theoretic intersections $D_j\cap D_k$ and $D_j\cap D_k \cap D_l$, and for intersections of cycles $D_jD_k$ and $D_jD_kD_l$. Over each $D_{jk}$, we get exceptional divisors $E_{jk,\alpha}$, $1\leq \alpha \leq s_{jk}$, where $s_{jk} = \ell(q_{jk},n)$ and whose components are over those of $D_{jk}$. From the local computations of the section above, we have
 \begin{equation}\label{eq:globpullD}
     h^*D_j = n\tilde{D}_j+ \sum_{D_kD_j \neq \emptyset}\sum_{\alpha = 1}^{s_{jk}}m_{jk,\alpha}E_{jk,\alpha} +\sum_{D_{kl}D_j \neq \emptyset}{F_{kl,j}},
 \end{equation}
 where $F_{kl,j}$ is a divisor whose components are the exceptional divisors over points in $D_{jkl}$.\\
 
 Explicitly, for any triple of positive integers numbers $j,k,l$, let $\rho_{kl}(j) \in \{1,2,3 \}$ the position of $j$ if we order the triple. For example $\rho_{23}(1)  = 1$, $\rho_{57}(6) = 2$ or $\rho_{54}(8) = 3$. Thus, if $F_{jkl,p}$ is the exceptional divisor over a point  $p\in D_{jkl}$, then
  $$F_{jkl} =  \sum_{p\in D_{jkl}} F_{jkl,p}$$
\begin{equation}\label{eq:Fklj}
 F_{kl,j} = \sum_{p\in D_{jkl}} v_{\rho_{kl}(j)}^{jkl}F_{jkl,p}.   
\end{equation}
 In terms of intersection theory, we have
 $$F_{jkl}^3 = \frac{n}{v^{jkl}_1v^{jkl}_2v^{jkl}_3}D_{jkl},\quad \tilde{D_j}\tilde{D_k} = 0.$$
\begin{remark}
    For the cyclic resolution constructed above, \Cref{sumnasym} gives
$$
V_{jkl}=0
$$
for every triple intersection. We keep the terms involving $V_{jkl}$
in the following formulas since they will also be useful for other
choices of the interior lattice points.
\end{remark}
 
 From \Cref{maincyclic} we have 
\begin{equation}\label{Qnumericaleq}
    K_{X_n} \sim_{\bQ} h^*K + \Delta,
\end{equation}
$$K := \left(K_Z+ \frac{n-1}{n}\sum_{j=1}^r D_j\right),\quad \Delta = \sum_{j<k} E_{jk} +\sum_{j<k<l}V_{jkl}F_{jkl},$$
$$E_{jk} = \sum_{\alpha=1}^{s_{jk}}N_{jk,\alpha}E_{jk,\alpha}\quad N_{jk,\alpha} = \frac{m_{jk,\alpha} + n_{jk,\alpha} -n}{n},\quad V_{jkl} = \frac{v_1^{jkl}+v_2^{jkl}+v_3^{jkl}-n}{n}.$$

Now we describe how the intersection theory on $X_n$ behaves under pullbacks of divisors of $Z$. In what follows, we set $E_{jk,0} = \tilde{D}_j$ and $E_{jk,s_{jk}+1} = \tilde{D}_k$.
 
\begin{prop}\label{pullbackextwall}
Let $G,G'$ any divisors on $Z$, then $h^*GF_{kl,j} = 0$ as $2$-cycle, for any $1\leq \alpha \leq s_{jk}$,
$$h^*Gh^*G'E_{jk,\alpha} = 0,$$
 $$h^*GE_{jk,\alpha}^2 = -k_{jk,\alpha}G D_{jk},\quad h^*G E_{jk,\alpha}E_{jk,\alpha\pm 1} = G D_{jk}.$$
\end{prop}
\proof{ We use the projection formula repeatedly. The first one is given by the fact that $h_* F_{jk,l}$ has codimension $3$. Now for a $1\leq \alpha\leq s_{jk}$ we have $h_*E_{jk,\alpha}$ supported in codimension $2$, thus 
$$ h^*Gh^*G'E_{jk,\alpha}=GG'h_*E_{jk,\alpha}  = 0.$$
Finally, for any $\alpha$ we have
$$h^* D_jh^*GE_{jk,\alpha} = 0 = h^*G(m_{jk,\alpha-1}E_{jk,\alpha-1}E_{jk,\alpha} + m_{jk,\alpha}E_{jk,\alpha}^2 + m_{jk,\alpha+1}E_{jk,\alpha+1}E_{jk,\alpha}),$$
$$h^* D_kh^*GE_{jk,\alpha} = 0 = h^*G(n_{jk,\alpha-1}E_{jk,\alpha-1}E_{jk,\alpha} + n_{jk,\alpha}E_{jk,\alpha}^2 + n_{jk,\alpha+1}E_{jk,\alpha+1}E_{jk,\alpha}).$$
The recursive relations with $k_{jk,\alpha}$ give
$$\left[ \begin{array}{cc}
m_{jk,\alpha} & m_{jk,\alpha+1}\\
n_{jk,\alpha} & n_{jk,\alpha+1}\\
\end{array}\right]  \left[ \begin{array}{c}
h^*G (E_{jk,\alpha}^2+k_{jk,\alpha}E_{jk,\alpha-1}E_{jk,\alpha})\\
h^*G (E_{jk,\alpha+1}E_{jk,\alpha}-E_{jk,\alpha-1}E_{jk,\alpha})
\end{array} \right]  = \left[ \begin{array}{c}
0\\
0 
\end{array}\right].$$
From \Cref{HJident} we have the determinant $m_{jk,\alpha}n_{jk,\alpha+1} -  m_{jk,\alpha+1}n_{jk,\alpha} = n$, thus
$$h^*G(E_{jk,\alpha}^2+k_{jk,\alpha}E_{jk,\alpha-1}E_{jk,\alpha}) = 0$$
$$h^*G(E_{jk,\alpha+1}E_{jk,\alpha}-E_{jk,\alpha-1}E_{jk,\alpha}) = 0 .$$
In particular, we have
$$h^*GE_{jk,\alpha+1}E_{jk,\alpha} = h^*G\tilde{D}_jE_{jk,1} = GD_{jk},$$
and the result follows.
\qed
}

In what follows, for simplicity, let us denote $K_X = K_{X_n}$.

\begin{corollary}\label{pullGcanon}
For any divisor $G$ on $Z$ we have
$$h^*GE_{jk}K_X = -D_{jk}G\left((N_{jk,1}+N_{kj,1}) + \sum_{\alpha=1}^{s_{jk}}(k_{jk,\alpha}-2)\right).$$
\end{corollary}
\proof{
Since $h^*G\, F$ vanishes at top-dimensional intersections, and $E_{jk,\alpha_1}E_{jk,\alpha_2} = 0$ for $|\alpha_1-\alpha_2| > 1$, we have
\begin{align*}
h^*GE_{jk}K_X &=h^*G(E_{jk})^2 = h^*G \sum_{\alpha=1}^{s_{jk}} N_{jk,\alpha}^2E_{jk,\alpha}^2 +2 N_{jk,\alpha}N_{jk,\alpha+1}E_{jk,\alpha}E_{jk,\alpha+1}   \\
&= D_{jk}G\sum_{\alpha=1}^{s_{jk}} -k_{jk,\alpha}N_{jk,\alpha}^2 +2 N_{jk,\alpha}N_{jk,\alpha+1}\\
&= D_{jk}G\sum_{\alpha=1}^{s_{jk}} N_{jk,\alpha}N_{jk,\alpha+1} - N_{jk,\alpha}N_{jk,\alpha-1} + N_{jk,\alpha}(k_{jk,\alpha}-2)\\
&= D_{jk}G\sum_{\alpha=1}^{s_{jk}}N_{jk,\alpha}(k_{jk,\alpha}-2)\\
&=-D_{jk}G \left((N_{jk,1}+N_{kj,1}) + \sum_{\alpha=1}^{s_{jk}}(k_{jk,\alpha}-2)\right).
\end{align*}
where the last identity is by telescoping sum argument.\\
\qed
}

\subsection{Asymptoticity of $K_{X_n}^3$.}\label{asympofK} Let us introduce the following notation 
$$|D|_{jk} := |D_j|_{k}+|D_k|_{j},$$
where $|D_j|_{k}$ satisfy
$$\sum_{D_{ll'}D_j\neq 0} F_{ll',j} E_{jk,\alpha}K_X = |D_j|_{k}(k_{jk,\alpha}-2).$$
\begin{lemma}\label{remarkasymp}
 We have
$$|D_j|_{k} = \sum_l \dfrac{v^{jkl}_{\rho_{kl}(j)}}{{v^{jkl}_{\rho_{kj}(l)}}}D_{jkl}.$$
\end{lemma}
\proof{Using equations from (\ref{eq:Kcurvasjk}), observe that $|D_j|_{k}$ depends on slopes of weights $v^{jkl}_1,v^{jkl}_2,v^{jkl}_3$ of the lattice points $v^{jkl}$. Using (\ref{eq:Fklj}), we get the result.\\
\qed
}

We need this to compute the intersection of $K_X$ with the external walls of the local toric resolution. Recursively we denote,
$$x_{jk,\alpha} = K_XE_{jk,\alpha-1}E_{jk,\alpha},\quad 1\leq \alpha \leq s_{jk}+1$$ 
$$y_{jk,\alpha} = K_XE_{jk,\alpha}^2,\quad 1\leq \alpha \leq s_{jk}.$$
Thus, we can write
\begin{equation}\label{eq:pullEk}
h^*D_j E_{jk,\alpha} K_{X} = m_{jk,\alpha-1}x_{jk,\alpha} + m_{jk,\alpha}y_{jk,\alpha} + m_{jk,\alpha+1}x_{jk,\alpha+1} + |D_j|_k(k_{jk,\alpha}-2)
\end{equation}
Using the $\bQ$-numerical equivalence of $K_X$ in \ref{maincyclic}, we compute
$$x_{jk,1} = h^*D_k\tilde{D_j}K_{X} = D_{jk}\left(K + \sum_{l\neq j} N_{jl,1}D_l \right), $$
$$x_{jk,s+1} = h^*D_j\tilde{D_k}K_{X} = D_{jk}\left(K + \sum_{l\neq k} N_{kl,1}D_l \right).$$

\begin{prop}\label{kcurvasglob}
We have
$$ x_{jk,\alpha} = x_{jk,1} + \frac{1}{n}(m_{jk,\alpha}^*(D_{jk}D_k-|D_k|_{j})  - n_{jk,\alpha}^*(D_{jk}D_j-|D_j|_{k}) ),$$
$$ y_{jk,\alpha} =  - k_{jk,\alpha}x_{jk,\alpha} + \frac{(k_{jk,\alpha} - 2)}{n}(n_{jk,\alpha+1}(D_{jk}D_j-|D_j|_{k}) - m_{jk,\alpha+1}(D_{jk}D_k-|D_k|_{j})).$$
Where $m_{jk,\alpha}^*=m_{jk,\alpha}-m_{jk,\alpha-1} - m_{jk,1}+m_{jk,0}$ and analogous for $n_{jk,\alpha}^*$.
\end{prop}
\proof{
Using the recursion given by the $k_{jk,\alpha}'s$, and formulas for $h^*D_j\, E_{jk,\alpha}K_X$  and $h^*D_k\, E_{jk,\alpha}K_X$ of (\ref{eq:pullEk}), we have
$$ \left[ \begin{array}{c}
(D_j^2D_k-|D_j|_{k})(k_{jk,\alpha}-2) \\
(D_jD_k^2-|D_k|_{j})(k_{jk,\alpha}-2)
\end{array}\right]  =  \left[ \begin{array}{cc}
m_{jk,\alpha} & m_{jk,\alpha+1}\\
n_{jk,\alpha} & n_{jk,\alpha+1}
\end{array}\right] \left[ \begin{array}{c}
k_{jk,\alpha}x_{jk,\alpha} + y_{jk,\alpha}\\
x_{jk,\alpha+1} - x_{jk,\alpha}
\end{array}\right]  $$
The determinant $m_{jk,\alpha}n_{jk,\alpha+1} -m_{jk,\alpha+1}n_{jk,\alpha} = n$, implies second relation for $y_{jk,\alpha}$, and
$$x_{jk,\alpha+1} = x_{jk,\alpha} + \frac{(k_{jk,\alpha} - 2)}{n}(m_{jk,\alpha}(D_{jk}D_k-|D_k|_{j})  - n_{jk,\alpha}(D_{jk}D_j-|D_j|_{k}) )$$
The recurrence for $x_{jk,\alpha}$ with a telescopic sum argument gives the result.\\
\qed
}

In the rest, we will denote
$$(k -2)_{jk} := \sum_{\alpha = 1}^{s_{jk}} (k_{jk,\alpha}-2).$$

\begin{prop}\label{computingK3}
    We have,
    \begin{align*}
K_X^3 = &nK^3  - 2\sum_{j<k} \left( D_{jk}K + \sum_{\substack{l\neq j,k\\ D_{jkl}\neq 0}} \frac{V_{jkl}}{v_{\rho_{jk}(l)}^{jkl}}D_{jkl}\right)(N_{jk,1}+N_{kj,1} + (k -2)_{jk})\\
&+ \sum_{j<k<l} \frac{nV_{jkl}^3}{v_1^{jkl}v_2^{jkl}v_3^{jkl}}D_{jkl} +\sum_{j<k}\sum_{\alpha=1}^{s_{jk}} \frac{D_{jk}(D_j+D_k)-|D|_{jk}}{n}N_{jk,\alpha}(k_{jk,\alpha}-2)\\
&-\sum_{j<k}\sum_{\alpha=1}^{s_{jk}}  N_{jk,\alpha}(x_{jk,\alpha} + y_{jk,\alpha}+x_{jk,\alpha+1}).
\end{align*}
\end{prop}
\proof{ We will compute $K_{X}^3$ using the numerical equivalence of (\ref{Qnumericaleq}). Squaring we get
\begin{align*}
K_{X}^2 \sim_{\bQ} h^*K^2 + \left( \sum_{j<k}E_{jk} \right)^2 + \sum_{j<k<l}V_{jkl}^2 F_{jkl}^2 + 2\sum_{j<k}E_{jk}\left(h^*K +\sum_{j<k<l}V_{jkl}F_{jkl}\right).    
\end{align*}
We have explicitly
$$(h^*K)^2K_X = (h^*K)^3 = nK^3.$$
$$E_{jk}F_{jkl}K_X = \frac{D_{jkl}}{v_{\rho_{jk}(l)}^{jkl}}\sum_{\alpha=1}^{s_{jk}}N_{jk,\alpha}(k_{jk,\alpha}-2) = -\frac{D_{jkl}}{v_{\rho_{jk}(l)}^{jkl}} \left((N_{jk,1}+N_{kj,1}) + (k -2)_{jk} \right).$$
Using \Cref{pullGcanon} for $G= K$, we get
\begin{align*}
K_X^3 = &nK^3  - 2\sum_{j<k} \left( D_{jk}K +  \sum_{\substack{l\neq j,k\\ D_{jkl}\neq 0}}  \frac{V_{jkl}}{v_{\rho_{jk}(l)}^{jkl}}D_{jkl}\right)(N_{jk,1}+N_{kj,1} + (k -2)_{jk})\\
&+ \sum_{j<k<l} \frac{nV_{jkl}^3}{v_1^{jkl}v_2^{jkl}v_3^{jkl}}D_{jkl} +K_X\left( \sum_{j<k}E_{jk} \right)^2.
\end{align*}
Just rest to compute
$$K_X\left( \sum_{j<k}E_{jk} \right)^2  = \sum_{j<k}  \sum_{\alpha=1}^{s_{jk}} N_{jk,\alpha}(N_{jk,\alpha-1}x_{jk,\alpha} + N_{jk,\alpha}y_{jk,\alpha} +N_{jk,\alpha+1}x_{jk,\alpha+1}).$$
From \Cref{pullGcanon}, we have
\begin{align*}
\frac{D_{jk}(D_j+D_k) (k_{jk,\alpha}-2)}{n} & =  \frac{h^*(D_j+D_k)E_{jk,\alpha}K_X}{n}\\
&= N_{jk,\alpha-1}x_{jk,\alpha} + N_{jk,\alpha}y_{jk,\alpha} +N_{jk,\alpha+1}x_{jk,\alpha+1}\\
 & + (x_{jk,\alpha} + y_{jk,\alpha}+x_{jk,\alpha+1}) +\frac{(k_{jk,\alpha}-2)}{n}|D|_{jk},
\end{align*}
So, we have explicitly
\begin{align*}
K_X\left( \sum_{j<k}E_{jk} \right)^2 &= \sum_{j<k}\sum_{\alpha=1}^{s_{jk}} \frac{D_{jk}(D_j+D_k)-|D|_{jk}}{n}N_{jk,\alpha}(k_{jk,\alpha}-2)\\
&-\sum_{j<k}\sum_{\alpha=1}^{s_{jk}}  N_{jk,\alpha}(x_{jk,\alpha} + y_{jk,\alpha}+x_{jk,\alpha+1}).
\end{align*}
\qed
}

\begin{theorem}\label{assymK}
If $\{D_1,\ldots,D_r \}$ is an asymptotic arrangement, then
$$\frac{K_{X_n}^3}{n} \to  -\logc_1^3(Z,D),$$ 
for prime numbers $n\to \infty$.
\end{theorem}
\proof{
The first term of the sum contains $\sum_{\alpha} N_{jk,\alpha}(k_{jk}-2)$, which is asymptotic respect to $n$ by previous discussion (\Cref{asympres2}). Thus, we just have to prove asymptoticity for
\begin{equation}\label{eq:firstred}
  \sum_{j<k}\sum_{\alpha=1}^{s_{jk}} N_{jk,\alpha}(x_{jk,\alpha} + y_{jk,\alpha}+x_{jk,\alpha+1})  
\end{equation}
Proceeding as above, it is not difficult to show the following identity,
\begin{align*}
\sum_{\alpha=1}^{s_{jk}} \frac{D_{jk}(D_j+D_k)-|D|_{jk}}{n}(k_{jk,\alpha}-2)&= \sum_{\alpha=1}^{s_{jk}} (N_{jk,\alpha}+1) (x_{jk,\alpha} + y_{jk,\alpha}+x_{jk,\alpha+1}) \\ 
& - N_{jk,1}x_{jk,1} - N_{jk,s}x_{jk,s+1}.
\end{align*}
So, the asymptoticity of (\ref{eq:firstred}) depends only on the asymptoticity of
\begin{equation}\label{eq:secondred}
 \sum_{j<k} \sum_{\alpha=1}^{s_{jk}}(x_{jk,\alpha} + y_{jk,\alpha}+x_{jk,\alpha+1}).
\end{equation}
By \Cref{kcurvasglob}, $x_{jk,\alpha} + y_{jk,\alpha}+x_{jk,\alpha+1}$ equals to
$$
 \frac{k_{jk,\alpha}-2}{n}(m_{jk,\alpha}^{**}(D_{jk}D_k-|D_k|_j)  - n_{jk,\alpha}^{**}(D_{jk}D_j-|D_j|_k) - nx_{jk,1} ),$$
where $m_{jk,\alpha}^{**} = m_{jk,\alpha-1} -m_{jk,\alpha+1} - m_{jk,0} + m_{jk,1}$, and analogous for $n_{jk,\alpha}^{**}$. Observe that these terms are asymptotically as $O(n)$. On the other hand, the terms $(D_{jk}D_j-|D_j|_k)$ and $(D_{jk}D_k-|D_k|_j)$ asymptotically depend only on the slopes of coordinates of the chosen lattice points $v^{jkl}$ on each intersection $D_{jkl}$. By \Cref{sumnasym}, we can choose lattice points with slopes asymptotically bounded by $C = 5$ as $n$ grows, with $K_X$ having $V_{jkl} = 0$ for all $j<k<l$. So, we have
$$|D|_{jk} \leq 10\sum_{l} D_{jkl}.$$
Therefore, \Cref{kcurvasglob} gives
\[
x_{jk,\alpha}+y_{jk,\alpha}+x_{jk,\alpha+1}
=O(k_{jk,\alpha}-2).
\]
Since $|N_{jk,\alpha}|\leq 1$ and
$\sum_{\alpha=1}^{s_{jk}}(k_{jk,\alpha}-2)=O(\sqrt n)$ by
\Cref{asymptoticOk}, the sum in (\ref{eq:firstred}) is $O(\sqrt n)$. Thus, as $n$ grows, all the terms in $K_X^3/n$  vanish except $K^3 \to -\logc_1^3(Z,D)$.\\
\qed}

\subsection{Asymptoticity of $e(X_n)$.}\label{asympofE}
The topological characteristic can be computed from the topology of $(Z,D)$ and the exceptional divisors $E_{jk,\alpha},F_{jkl}$. The divisors $F_{jkl} = \sum_{p\in D_{jkl}} F_{jkl,p}$ where $F_{jkl,p}$ is the corresponding exceptional divisor over a $p \in D_{jkl}$. Thus, $e(F_{jkl}) = D_{jkl}e(F_{jkl,p})$. In the toric picture (\Cref{toricreso}) of 
 $X_n$ over $p$, let $v = v^{jkl}$ the ray generator defining $F_{jkl,p}$. 
 \begin{lemma}\label{topofF}
We have $e(F_{jkl,p}) =  s_{jk}+s_{jl} +s_{kl} + 3$.
 \end{lemma}
 \proof{
It is well-known that $F_{jkl,p}$ is the toric variety associated to the star-fan $\mbox{Star}(C(v))\subset (N_v)_{\bR}$, i.e., the induced fan by the lattice quotient $N_{v} = N/vN$. In this case, the $2$-cones of  $\mbox{Star}(C(v))$ are the induced by each $C(e_{jk,\alpha},e_{jk,\alpha+1})$. Since $e(F_{jkl,p})$ is the sum of its top-dimensional cones \cite[Thm. 12.3.9]{coxlittleschenk2011}, we have the result.\\
\qed
 }

\begin{lemma}
If $X$ is a complex algebraic variety, and $A\subset X$ is a subvariety such that $X\setminus A$ is smooth, then $e(X) = e(X\setminus A) + e(A)$. More generally, the topological Euler characteristic is additive with
respect to finite constructible stratifications.
\end{lemma}
\proof{See \cite[p. 141]{fulton1993introduction}.\\ \qed}

\begin{prop}
    We have
$$
e(X_n)
=
n\bigl(e(Z)-e(D_{\mathrm{red}})\bigr)
+
e(D_{\mathrm{red}})
+
\sum_{j<k}s_{jk}\,e(D_{jk})
+
2\sum_{j<k<l}D_{jkl},
$$
where $s_{jk}=\ell(q_{jk},n)$ is the length of the Hirzebruch--Jung
resolution over $D_{jk}$.
\end{prop}
\proof{
For each $j<k$, set
\[
D_{jk}^{\circ}
:=
D_{jk}\setminus \bigcup_{l\neq j,k}D_{jkl},
\]
and let $D^{\circ}$ be the locus of $D_{red}$ contained in exactly one
irreducible component. Then $Z$ admits the disjoint stratification
$$
Z
=
(Z\setminus D_{\mathrm{red}})
\sqcup D^{\circ}
\sqcup \bigsqcup_{j<k}D_{jk}^{\circ}
\sqcup \bigsqcup_{j<k<l}D_{jkl}.
$$
We compute the Euler characteristic of the fibers of
\[
h:X_n\longrightarrow Z
\]
over each stratum. 
\begin{enumerate}
    \item Outside $D_{red}$, the morphism $h$ is unramified of degree $n$, hence
$$
e\bigl(h^{-1}(Z\setminus D_{\mathrm{red}})\bigr)
=
n\,e(Z\setminus D_{red}).
$$

\item Over $D^{\circ}$ we have a homeomorphism, so $$e(h^{-1}(D^{\circ})) = e(D^{\circ})$$. 

\item Over
$D_{jk}^{\circ}$ the fiber is a chain of $s_{jk}$ rational curves. Therefore
$$
e\bigl(h^{-1}(x)\bigr)
=
2s_{jk}-(s_{jk}-1)
=
s_{jk}+1,
\qquad x\in D_{jk}^{\circ}.
$$

\item Finally, if $p\in D_{jkl}$, the fiber is the toric exceptional divisor
$F_{jkl,p}$. By \Cref{topofF},
$$
e(F_{jkl,p})
=
s_{jk}+s_{jl}+s_{kl}+3.
$$
\end{enumerate}
By additivity of the topological Euler characteristic, we obtain
$$
e(X_n)
=
n\,e(Z\setminus D_{\mathrm{red}})
+e(D^{\circ}) 
+
\sum_{j<k}(s_{jk}+1)e(D_{jk}^{\circ}) 
+
\sum_{j<k<l}
(s_{jk}+s_{jl}+s_{kl}+3)D_{jkl}.
$$
On the other hand,
$$
e(D_{\mathrm{red}})
=
e(D^{\circ})
+
\sum_{j<k}e(D_{jk}^{\circ})
+
\sum_{j<k<l}D_{jkl},
$$
and
$$
e(D_{jk}^{\circ})
=
e(D_{jk})
-
\sum_{l\neq j,k}D_{jkl}.
$$
Therefore,
$$
e(X_n)
=
n\bigl(e(Z)-e(D_{\mathrm{red}})\bigr)
+e(D_{\mathrm{red}})
+\sum_{j<k}s_{jk}e(D_{jk}^{\circ}) 
+
\sum_{j<k<l}
(s_{jk}+s_{jl}+s_{kl}+2)D_{jkl}.
$$
Substituting the expression for $e(D_{jk}^{\circ})$, every triple
intersection $D_{jkl}$ contributes
$$
-(s_{jk}+s_{jl}+s_{kl})D_{jkl}
+
(s_{jk}+s_{jl}+s_{kl}+2)D_{jkl},
$$
so the terms depending on the lengths cancel. This proves the formula.\\
\qed
}

As a corollary of the previous identity for $e(X_n)$, and
by the previous discussion in \Cref{asympres2}, the lengths $s_{jk}/n$ are asymptotically zero as $n$ grows by \Cref{girsttheor}. Thus, we get the following.

\begin{corollary}\label{asympE}
If $\{D_1,\ldots,D_r \}$ is an asymptotic arrangement, then for prime numbers $n\to \infty$
$$\frac{e(X_n)}{n} \longrightarrow \logc_3(Z,D).$$ 
\end{corollary}
\qed

\bibliographystyle{alpha}
\bibliography{main.bib}

\end{document}